\documentclass[12pt]{amsart}

\usepackage[utf8]{inputenc}
\usepackage[T1]{fontenc}
\usepackage{amssymb, amsthm, amsmath}
\usepackage{enumitem}
\usepackage{hyperref}
\usepackage{comment}
\usepackage{centernot}

\usepackage{tikz}
\usetikzlibrary{arrows.meta,fit}
\newcommand{\core}{\text{\textcircled{\tiny\#}}}

\numberwithin{equation}{section}

\newcommand{\nimplies}{\;\centernot\Longrightarrow\;}

\newcommand{\PP}{\mathsf{P}}
\newcommand{\QQ}{\mathsf{Q}}

\def\N{{\mathcal{N}}}
\def\R{{\mathcal{R}}}

\def\*<{\,*\!<}
\def\<*{< \!*\,}

\def\G{\mathcal{G}}
\def\r{\operatorname{rank}}

\newtheorem{theorem}{Theorem}[section]
\newtheorem{lemma}[theorem]{Lemma}
\newtheorem{proposition}[theorem]{Proposition}
\newtheorem{corollary}[theorem]{Corollary}
\theoremstyle{definition}
\newtheorem{definition}[theorem]{Definition}
\newtheorem{example}[theorem]{Example}
\newtheorem{remark}[theorem]{Remark}

\begin{document}

\title[SQUARE ORDER AND THE PLUS ORDER]{A square analogue of the diamond order and a unified approach to the plus order}

\author{Dragan S. Raki\'c}
\address{Faculty of Mechanical Engineering, University of Ni\v s, Aleksandra Medvedeva 14, 18000 Ni\v s, Serbia}
\email{rakic.dragan@gmail.com}

\begin{abstract}
We introduce the square order $\le^\square$ on group invertible matrices, defined by $A^3=ABA$ and $A\le^s B$, where $\le^s$ is the space pre-order. It is the group inverse analogue of the diamond order: $A\le^\square B$ if and only if $A^\#\le^- B^\#$.
Using its block representation, we characterize the reverse order law $(AB)^\#=B^\#A^\#$ for $A\le^\square B$.
In contrast to the Moore--Penrose case, the conditions $A\le^- B$ and $(B-A)^\#=B^\#-A^\#$ do not characterize the sharp order; instead, they lead to new partial orders between the sharp and square orders.
We also develop a framework where the plus order replaces the minus order in Mitra's unified theory. For a map $\mathcal G$ assigning to each matrix $A$ a subset of its inner inverses, we define $A\le^{\mathcal G}_{+}B$ if $A\le^s B$ and $A=PBQ$ for some pair $(P,Q)=(AA^-,A^-A)$, $A^-\in\mathcal G(A)$, whereas Mitra's relation requires $A=PB=BQ$, and $A\preceq^{\mathcal G}_{+}B$ if $A\le^s B$ and $A^-\le^-B^-$ for some reflexive $A^-\in\mathcal G(A)$ and $B^-\in\mathcal G(B)$.
When $\mathcal G(A)$ consists of a single reflexive inverse, we characterize the transitivity of $\le^{\mathcal G}_{+}$ in terms of $\preceq^{\mathcal G}_{+}$, and in general we give a sufficient condition for both to be partial orders. For maps defined by the equations of the Moore--Penrose, group, core and dual core inverses, both relations coincide and only nine arise: the plus, diamond and square orders, two new orders based on core inverses, and four further relations. In Mitra's theory, they yield exactly the minus, star, sharp, core and dual core orders and the one-sided star and sharp orders.
\end{abstract}



\maketitle

\section{Introduction and preliminaries}

Partial orders on matrices defined by means of generalized inverses have
been studied extensively since the work of Hartwig~\cite{Ha},
Nambooripad~\cite{Nambooripad} and Mitra~\cite{Mitra_sharp_order}; see the monograph~\cite{knjiga}
and the references therein. The most fundamental of them is the minus
order, while the star, sharp and core orders are obtained by replacing an
arbitrary inner inverse in its definition by the Moore--Penrose, group and
core inverse, respectively. Mitra~\cite{Unified_theory} unified these orders by means
of g-maps. A different line of research started with the diamond order of
Baksalary and Hauke~\cite{Baksalary_Hauke_diamond}, which is characterized by
$A^\dagger\le^- B^\dagger$, and led to the plus order recently introduced
by Arias and Maestripieri~\cite{Arias_Maestripieri}.

In this paper we introduce the square order $\le^\square$ on the set of
group invertible matrices, defined by $A\le^s B$ and $A^3=ABA$. It is the
group inverse analogue of the diamond order: $A\le^\square B$ if and only
if $A^\#\le^- B^\#$ (Theorem~\ref{th_square_equiv_minus}). We obtain a block
representation of $\le^\square$ (Corollary~\ref{square_matrix_form}) and show that,
for $A\le^\square B$, the reverse order law $(AB)^\#=B^\#A^\#$ holds if and
only if $A\le^\# B$ (Theorem~\ref{obrni_okreni}). We also study additive
properties of the group inverse. In contrast to the Moore--Penrose case,
where $A\le^- B$ and $(B-A)^\dagger=B^\dagger-A^\dagger$ characterize the
star order~\cite{Baksalary_aditivity_MP,Fongi_Gonzalez,Hartwig_Styan}, the corresponding conditions for the group inverse do
not characterize the sharp order; instead, they lead to new partial orders
$\le^\boxplus_\ell$, $\le^\boxplus_r$ and $\le^\boxplus$
(Corollary~\ref{additivity_group_cor}).

In Section~3 we develop a framework in which the plus order plays the
role of the minus order in Mitra's theory.
For a g-map $\G$, which assigns to each matrix $A$ a (possibly empty) subset $\G(A)$ of $A\{1\}$, we introduce two relations: the \emph{g-plus order relation}, $A\le^{\G}_{+}B$ if $A\le^s B$ and $A=AA^-BA^-A$ for some $A^-\in\G(A)$, and the \emph{strong g-plus order relation}, $A\preceq^{\G}_{+}B$ if $A\le^s B$ and $A^-\le^-B^-$ for some reflexive $A^-\in\G(A)$ and $B^-\in\G(B)$.
The plus, diamond and square orders are special cases of
both relations. We describe the maximal elements and, when each $\G(A)$
consists of a single reflexive inverse, we prove that $\le^{\G}_{+}$ is
transitive if and only if $A\le^{\G}_{+}B$ implies $A\preceq^{\G}_{+}B$
whenever $B$ is not maximal (Theorem~\ref{th_<=preceq}); this is an analogue of
a result of Mitra~\cite[Theorem~2.6]{Unified_theory}. For general g-maps we
give a sufficient condition under which the two relations coincide and
are partial orders (Theorem~\ref{th_dovoljan_uslov}).

We then consider the g-maps $\G(A)=A\{1,S\}$ defined by a fixed set $S$
of the equations defining the Moore--Penrose, group, core and dual core
inverses. We show that all of them satisfy this sufficient condition
and that, although there are $24$ such g-maps, they induce only nine
relations, four of them if $m\neq n$
(Theorem~\ref{teorema_klasifikacija}, Table~\ref{tabela_relacija}).
Besides the plus, diamond and square orders, these include two new partial
orders $\le^\bigcirc$ and $\le_\bigcirc$ related to the core and dual core
inverses, and four further relations. The key observation is that
$\le^{\G}_{+}$ depends only on the pairs of idempotents $(AA^-,A^-A)$,
$A^-\in\G(A)$: $A\le^{\G}_{+}B$ if and only if $A\le^s B$ and $A=PBQ$ for
such a pair $(P,Q)$, whereas Mitra's relation is given by $A=PB=BQ$.
The same approach applied to Mitra's relations yields the star, core,
dual core, sharp, one-sided star, one-sided sharp and minus orders
(Table~\ref{tabela_Mitra}), and shows that these are the only orders induced by the g-maps $A\{1,S\}$ in
Mitra's sense. We conclude by comparing the two theories.

We now fix the notation and recall the necessary definitions.

Let $M_{m\times n}=\mathbb{C}^{m\times n}$ denote the set of all complex matrices of size $m\times n$ and let $A \in M_{m\times n}$.
We denote by $\R(A)$ and $\N(A)$ the range and null space of $A$, respectively.
Also, $A^{*}$ is the conjugate transpose of $A$.

For $A\in M_{m\times n}$, a matrix $G \in M_{n\times m}$ is called a \emph{generalized inverse} (or g-inverse or $\{1\}$-inverse) of $A$ if it satisfies
$AGA = A.$
In general, a generalized inverse of $A$ is not unique; the set of all generalized inverses of $A$ is denoted by $A\{1\}$.
If in addition $GAG=G$, then $G$ is a \emph{reflexive g-inverse} of $A$.
For $A\in M_{m\times n}$ consider the equations
\begin{align*}
(1)\;\;&AX A=A, &
(2)\;\;&XAX=X,\\
(3)\;\;&(AX)^*=AX, &
(4)\;\;&(XA)^*=XA.
\end{align*}
For every $A\in M_{m\times n}$ there exists a unique matrix
$A^\dag\in M_{n\times m}$, called the \emph{Moore--Penrose inverse} of $A$, that
satisfies equations (1)--(4).
For brevity, let $M_n=M_{n\times n}$ denote the set of all complex square matrices of order $n$.
The \emph{index} of a matrix $A\in M_n$, denoted by $\operatorname{ind}(A)$, is defined as the smallest nonnegative integer $k$ such that
$$
\operatorname{rank}(A^{k}) = \operatorname{rank}(A^{k+1}).
$$
A matrix $A$ is nonsingular (invertible) if and only if $\operatorname{ind}(A)=0$.
When $\operatorname{ind}(A)\le 1$, there exists a unique matrix
$A^\#\in M_n$, called the \emph{group inverse} of $A$, that satisfies (1), (2) and
$$\text{(5) } AX=XA.$$
We denote the set of all group invertible matrices by $M_n^\#$.
A matrix $A\in M_n$ is called \emph{EP} if $\R(A)=\R(A^*)$ (equivalently, $\mathcal N(A)=\mathcal N(A^*)$); in this case $\operatorname{ind}(A)\le 1$. It is well known that $A$ is EP if and only if $A\in M_n^\#$ and $A^\dag=A^\#$.
For $A\in M_{m\times n}$ and $i,j,\dots,k\in\{1,2,\dots\}$, we denote by
$A\{i,j,\dots,k\}$ the set of all matrices $X\in M_{n\times m}$ satisfying
equations $(i),(j),\dots,(k)$.

Generalized inverses provide the algebraic foundation for several partial orders on matrices, the most fundamental of which is the minus order.
For $A, B \in M_{m\times n}$, we say that $A$ is below $B$ with respect to the \emph{minus partial order}, denoted $A \le^{-} B$, if there exists a generalized inverse $A^{-} \in A\{1\}$ such that
$$A^{-} A = A^{-} B \qquad \text{and} \qquad A A^{-} = B A^{-}.$$
We will use the following important characterizations of the minus partial order.
\begin{theorem}\cite{knjiga}\label{minus_karakterizacije}
Let $A,B\in M_{m\times n}$. Then the following assertions are equivalent:
\begin{enumerate}[label=(\roman*)]
\item $A\le^- B$;
\item $\operatorname{rank}(B-A) = \operatorname{rank}(B)-\operatorname{rank}(A)$;
\item $\R(A)\oplus \R(B-A)=\R(B)$,
where $\oplus$ denotes the direct sum of subspaces;
\item $B\{1\}\subseteq A\{1\}$;
\item $A=AB^-A=AB^-B=BB^-A$ for every $B^-\in B\{1\}$;
\item $A=AB^-A=AB^-B=BB^-A$ for some $B^-\in B\{1\}$.
\end{enumerate}
\end{theorem}
The minus order was introduced in 1980 by Hartwig \cite{Ha} and, independently, by Nambooripad \cite{Nambooripad}, and it represents one of the earliest and most extensively studied partial orders on the set of matrices.
Several important orders refine the minus order, including the star order, the sharp order, and their one-sided (left and right) counterparts.

The \emph{star partial order} on $M_{m\times n}$, denoted by $\le^*$, is defined by
$$
A\le^* B \iff AA^\dag=BA^\dag \text{ and } A^\dag A=A^\dag B.
$$
The star partial order was introduced by Drazin \cite{Drazin_natural_structures} in 1978.

On the set $M_n^\#$, Mitra \cite{Mitra_sharp_order} introduced in 1987 the so-called \emph{sharp partial order}, defined as follows: for $A, B \in M_n^\#$, we say that $A \le^{\#} B$ if and only if
$$
A^{\#} A = A^{\#} B \qquad \text{and} \qquad A A^{\#} = B A^{\#}.
$$

Later, as a natural weakening of the sharp order, the so-called \emph{left-sharp} and \emph{right-sharp} partial orders were introduced by Mitra et al. in \cite{knjiga}. For $A,B\in M_n^\#$,
\begin{align*}
A \le^{\#}_{\ell} B & \iff A^{\#} A = A^{\#} B \;\text{ and }\; \R(A)\subseteq\R(B), \\
A \le^{\#}_{r} B & \iff A A^{\#} = B A^{\#} \;\text{ and }\; \R(A^*)\subseteq\R(B^*).
\end{align*}
It is clear that $A \le^{\#} B$ if and only if $A \le^{\#}_{\ell} B$ and $A \le^{\#}_{r} B$.
It is worth noting that, if $A, B \in M_n^{\#}$ then
$$
A^{\#} A = A^{\#} B \iff A^2 = AB,
$$
and
$$
A A^{\#} = B A^{\#} \iff A^2 = BA.
$$
Also, for $A, B\in M_{m\times n}$,
$$A^{\dag} A = A^{\dag} B \iff A^*A = A^*B, $$
and
$$A A^{\dag} = B A^{\dag} \iff AA^* = BA^*.$$

The fundamental relation on $M_{m\times n}$, closely related to the orders considered in this paper, is the \emph{space pre-order}. For $A, B \in M_{m\times n}$, it is defined by
$$
A \le^{s} B \iff \R(A) \subseteq \R(B) \;\text{ and }\; \R(A^{*}) \subseteq \R(B^{*}).
$$
This relation is reflexive and transitive, but not antisymmetric in general, and therefore defines only a pre-order on $M_{m\times n}$. It plays an important role, since all relations discussed in this paper imply the space pre-order.
Note that $\R(A) \subseteq \R(B)$ if and only if $A=BX$ for some $X\in M_n$ if and only if $A=BB^-A$ for every $B^-\in B\{1\}$. Also, $\R(A^{*}) \subseteq \R(B^{*})$ if and only if $\N(B) \subseteq \N(A)$ if and only if $A=YB$ for some $Y\in M_m$ if and only if $A=AB^-B$ for every $B^-\in B\{1\}$.

Denote by $E(M_n)$ and  $E^*(M_n)$ the sets of idempotent and self-adjoint idempotent matrices in $M_n$, respectively, i.e., $$E(M_n)=\{A\in M_n:A^2 = A\} \text{ and } E^*(M_n)=\{A\in M_n:A^2=A=A^*\}.$$

\section{Square partial order}

Before introducing a new relation, we recall the definition of the \emph{diamond partial order} on $M_{m\times n}$, denoted by $\le^\diamondsuit$. This partial order was introduced by Baksalary and Hauke \cite{Baksalary_Hauke_diamond} in 1990 and is defined by
$$
A\le^\diamondsuit B \iff A\le^s B \text{ and } AA^*A=AB^*A.
$$
Moreover, the following characterization holds:
\begin{equation}\label{diamond_kao_preceq}
A\le^\diamondsuit B \iff A^\dag\le^- B^\dag.
\end{equation}
The diamond partial order has been widely studied not only for matrices, but also in the setting of bounded linear operators on Hilbert spaces, $C^*$-algebras, and various classes of rings; see, e.g., \cite{Arias_Maestripieri,Baksalary_Hauke_diamond,BurgosMGMC,CirulisDiamond,FerreyraLevisMaharanaOrquera,FerreyraMalikCore,
Fongi_Gonzalez,HernandezLattanziThome_Diamond,LebtahiPatricioThome,MalikRuedaThome,MarovtMihelic,YangJi}.

\subsection{Definition and basic properties.}

Inspired by the diamond partial order, we introduce the following relation.
\begin{definition}\label{square_definition}
Let $A,B\in M_n$. We say that $A\le^\square B$ if $A,B\in M_n^\#$ and
$$A\le^s B \quad \text{ and } \quad A^3=ABA.$$
We call $\le^\square$ the \emph{square order}.
\end{definition}

The assumption $B\in M_n^{\#}$ in Definition \ref{square_definition} does not follow from the remaining requirements.
(To see this take $A=\operatorname{diag}(1,0,0)$ and $B=\begin{bmatrix}1&0&0\\0&0&1\\0&0&0\end{bmatrix}$).
We include this condition for two reasons. First, it makes $\le^\square$ a relation on the set $M_n^\#$, as is the case for the sharp order. Second, it guarantees that the group inverse $B^\#$ exists. We use $B^\#$ throughout the paper, for example in Theorem~\ref{th_square_equiv_minus}.

Since $\R(A^*)^{\perp}=\N(A)$, we have
$A\le^s B$ if and only if $\R(A)\subseteq\R(B)$ and $\N(B)\subseteq\N(A)$, which is equivalent to $A^*\le^s B^*$. Consequently, $A\le^{\square}B$ if and only if $A^*\le^\square B^*$. Moreover, $\le^\square$ is invariant under similarity, i.e. $A\le^\square B$ if and only if $TAT^{-1}\le^\square TBT^{-1}$ for every invertible matrix $T\in M_n$.

The following theorem characterizes the relation $\le^\square$ in terms of group inverses and the minus partial order.
\begin{theorem}\label{th_square_equiv_minus}
Let $A,B\in M_n^\#$. Then
$$A\le^\square B \iff A^\#\le^- B^\#.$$
\end{theorem}
\begin{proof}
Assume that $A\le^\square B$. Then $A^3=ABA$ and $A=BX=YB$, for some $X,Y\in M_n$. For arbitrary $G\in B^\#\{1\}$ we have
\begin{align*}
A^\# GA^\# &=(A^\#)^2AGA (A^\#)^2=(A^\#)^2 YBGBX (A^\#)^2 \\
&=(A^\#)^2 YB^2 B^\# G B^\# B^2 X (A^\#)^2=(A^\#)^2 YB^2 B^\# B^2 X (A^\#)^2 \\
&= (A^\#)^2 ABA(A^\#)^2 =(A^\#)^2 A^3 (A^\#)^2=A^\#,
\end{align*}
so $G\in A^\#\{1\}$. In view of Theorem~\ref{minus_karakterizacije}, we conclude that $A^\#\le^- B^\#$.

Conversely, assume that $A^\#\le^- B^\#$. Since $B\in B^\#\{1\}$, by Theorem~\ref{minus_karakterizacije}, we have that $A^\#BA^\#=A^\#$. Multiplying this equation on both sides by $A^2$, we obtain $ABA=A^3$. Since $A^\#\le^- B^\#$ implies $A^\#\le^s B^\#$, and since $\R(T)=\R(T^\#)$ and $\N(T)=\N(T^\#)$, for any $T\in M_n^\#$, we conclude that $A\le^s B$. Thus, $A\le^\square B$.
\end{proof}

\begin{corollary}\label{cor_square_je_pu}
The relation $\le^\square$ is a partial order on $M_n^\#$.
\end{corollary}
\begin{proof}
Follows by Theorem~\ref{th_square_equiv_minus} since the map $A\mapsto A^{\#}$ is an involution on $M_n^{\#}$ and $\le^{-}$ is a partial order.
\end{proof}

\begin{corollary}\label{square_preko_rangova}
Let $A,B\in M_n^{\#}$.
\begin{enumerate}[label=(\roman*)]
\item $A\le^{\square}B \iff A^{\#}=A^{\#}BA^{\#}=A^{\#}BB^{\#}=B^{\#}BA^{\#}$.
\item $A\le^{\square}B \iff \mathcal{R}(B)=\mathcal{R}(A)\oplus\mathcal{R}(B^{\#}-A^{\#})
      \iff \operatorname{rank}(B)=\operatorname{rank}(A)+\operatorname{rank}(B^{\#}-A^{\#})$.
\item $A\le^{-}B \iff A^{\#}\le^{\square}B^{\#}$.
\item  If $\operatorname{rank}(B-A)=\operatorname{rank}(B^\#-A^\#)$, then $A\le^\square B \iff A\le^-B$.
Conversely, if $A\le^-B$ and $A\le^\square B$, then $\operatorname{rank}(B-A)=\operatorname{rank}(B^\#-A^\#)$.
\item If $B^\#-A^\#\in M_n^\#$, then $A\le^\square B \iff (B^\#-A^\#)^\#\le^\square B$.
\end{enumerate}
\end{corollary}

\begin{proof}
(i)
Apply Theorem~\ref{th_square_equiv_minus} and Theorem~\ref{minus_karakterizacije}(vi) to $A^{\#}\le^{-}B^{\#}$, with $B\in B^{\#}\{1\}$.

(ii) Follows from Theorem~\ref{minus_karakterizacije}(ii),(iii), since
$\mathcal{R}(T)=\mathcal{R}(T^{\#})$ and $\operatorname{rank}(T)=\operatorname{rank}(T^{\#})$
for $T\in M_n^{\#}$.

(iii) Follows by Theorem~\ref{th_square_equiv_minus}.

(iv)
If $\operatorname{rank}(B-A)=\operatorname{rank}(B^{\#}-A^{\#})$, then, by (ii) and Theorem~\ref{minus_karakterizacije}(ii), both $A\le^{\square}B$ and $A\le^{-}B$ are equivalent to the same rank equality.
The converse follows since both $A\le^-B$ and $A\le^\square B$ give $\operatorname{rank}(B)-\operatorname{rank}(A)$ for the two ranks.

(v) By Theorem~\ref{th_square_equiv_minus},
$(B^\#-A^\#)^\#\le^\square B \iff B^\#-A^\#\le^- B^\# \iff A^\#\le^- B^\# \iff A\le^\square B$.
\end{proof}

Recall that $A\le^- B$ implies $B-A\le^- B$. By Corollary~\ref{square_preko_rangova}(v), the analogous statement for the square order is that $A\le^\square B$ implies $(B^\#-A^\#)^\#\le^\square B$ whenever $B^\#-A^\#\in M_n^\#$, while $B-A\le^\square B$ need not hold. Indeed, let $A=\operatorname{diag}(1,0)$ and $B=\begin{bmatrix}1&1\\1&0\end{bmatrix}$. Then $A\le^\square B$ and $B-A\in M_2^\#$, but $(B-A)^3=\begin{bmatrix}0&1\\1&0\end{bmatrix}\neq \begin{bmatrix}0&1\\1&1\end{bmatrix}=(B-A)B(B-A)$, so $B-A\le^\square B$ does not hold.

The \emph{natural partial order} $\le$ on $E(M_n)$ is defined in the following way. For $E,F\in E(M_n)$
$$E\le F \iff E=EF=FE.$$
For $A\in M_n^\#$, denote by $Q_A$ the projection onto $\R(A)$ along $\N(A)$. That is, $Q_A$ is the idempotent matrix
$$Q_A=AA^\#.$$

\begin{proposition}\label{square_sledi_plus}
Let $A,B\in M_n^\#$.
Then
$$Q_A\le Q_B \iff A\le^s B \iff A=Q_BAQ_B,$$
and the following statements are equivalent:
\begin{enumerate}[label=(\roman*)]
\item $A\le^\square B$;
\item $A\le^s B$ and $A=Q_ABQ_A$;
\item $A\le^s B$ and $A=QBQ$ for some $Q\in E(M_n)$;
\item $A\le^s B$ and $A^\#=A^\#BA^\#$.
\end{enumerate}
\end{proposition}
\begin{proof}
The proofs of the first statement and implications (i) $\Rightarrow$ (ii) $\Rightarrow$ (iii) and (ii) $\Leftrightarrow$ (iv) are straightforward. Assume that $A\le^s B$ and $A=QBQ$ for some $Q\in E(M_n)$. Then $A=QA=AQ$, so $ABA=AQBQA=A^3$. Thus, (iii) implies (i).
\end{proof}

\begin{proposition}\label{prop:basic}
Let $A,B,C\in M_n^\#$ and $E,F\in E(M_n)$.
\begin{enumerate}[label=(\roman*)]
\item $E\le^\square F \iff E\le^\#F \iff E\le^-F \iff E\le^s F \iff E\le F$.
\item $A\in E(M_n) \iff A\le^\square I_n \iff A\le^\square F$ for some $F\in E(M_n)$.
\item $A$ is a maximal element in $(M_n^\#,\le^\square)$ if and only if $A$ is invertible.
\item If $A\le^\square C$ and $B\le^ \square C$, then $A\le^ \square B \iff A\le^s B$.
\end{enumerate}
\end{proposition}

\begin{proof}
Item (i) is straightforward.

(ii) Since $A\le^s I_n$, we have $A\le^{\square}I_n$ if and only if $A^3=A^2$, which is equivalent to $A^2=A$ (multiply by $A^\#$). If $A\le^{\square}F$ with $F\in E(M_n)$, then $F\le^{\square}I_n$, so $A\le^{\square}I_n$ by transitivity (Corollary~\ref{cor_square_je_pu}). The converse follows with $F=I_n$.

(iii) If $A$ is invertible and $A\le^{\square}B$, then $A^3=ABA$ gives $A=B$.
Conversely, suppose that $A$ is not invertible. Then $B:=A+I_n-AA^{\#}\neq A$ is invertible, with
$B^{-1}=A^{\#}+I_n-AA^{\#}$, and $ABA=A^3$, so $A\le^{\square}B$.

(iv) Assume that $A\le^\square C$ and $B\le^\square C$. Then $A^3=ACA$ and $B^3=BCB$.
If $A\le^{s}B$, then $A=BX=YB$ for some $X,Y\in M_n$, and $ABA=YB^3X=YBCBX=ACA=A^3$, so $A\le^\square B$. The converse implication is trivial.
\end{proof}

\subsection{Block representation.}

In this subsection we describe $\le^\square$ by means of block representations with respect to a suitable decomposition of $\mathbb{C}^n$, where matrices are identified with the corresponding linear operators on $\mathbb{C}^n$. These representations are the main tool in the study of the reverse order law (Theorem~\ref{obrni_okreni}) and of additive properties of the group inverse (Theorem~\ref{additivity_group_th}).

\begin{lemma}\label{lema_space_order_decomposition}
Let $A,B\in M_n^\#$. Then $A\le^s B$ if and only if
\begin{align}
\R(B)&=\R(A)\oplus\left(\N(A)\cap\R(B)\right) \label{prva} \text{ and } \\
\N(A)&=\left(\N(A)\cap\R(B)\right)\oplus\N(B)\label{druga}.
\end{align}
Consequently,
\begin{equation}\label{dekomp_prostora}
\mathbb{C}^n=\R(A)\oplus\left(\N(A)\cap\R(B)\right)\oplus\N(B).
\end{equation}
Moreover, in this case,
\begin{equation}\label{treca}
\N(A)\cap\R(B)=\R(BB^\#-AA^\#).
\end{equation}
\end{lemma}
\begin{proof}
Since $A,B\in M_n^\#$, we know that
$$\mathbb{C}^n=\R(A)\oplus \N(A)=\R(B)\oplus \N(B).$$
Assume that $A\le^s B$, that is, $\R(A)\subseteq \R(B)$ and $\N(B)\subseteq \N(A)$.
We prove \eqref{prva}; the proof of \eqref{druga} is analogous.
It is easy to see that $\R(A)\cap\left(\N(A)\cap\R(B)\right)=\{0\}$ and $\R(A)\oplus\left(\N(A)\cap\R(B)\right)\subseteq\R(B)$.
To show that
$\R(B)\subseteq\R(A)\oplus\left(\N(A)\cap\R(B)\right)$,
let $x\in \R(B)$ and write $x=x_1+x_2$, where $x_1\in\R(A)$ and $x_2\in\N(A)$. Then $x_2=x-x_1\in \R(B)$, so $x=x_1+x_2\in\R(A)\oplus(\N(A)\cap\R(B))$.
Therefore, the decomposition~\eqref{prva} follows, and consequently so does \eqref{dekomp_prostora}.
Conversely, if \eqref{prva} and \eqref{druga} hold, then it is clear that $\R(A)\subseteq\R(B)$ and $\N(B)\subseteq \N(A)$.

It remains to prove that \eqref{treca} holds whenever $A\le^s B$. Let $x\in \N(A)\cap\R(B)$, that is, $Ax=0$ and $x=Bz$ for some $z\in\mathbb{C}^n$. Then
$(BB^\#-AA^\#)x=BB^\# Bz-A^\#Ax=Bz=x$, so $x\in \R(BB^\#-AA^\#)$. Conversely, since $A\le^s B$, we have $A=BX=YB$ for some $X,Y\in M_n$. Moreover, $BB^\#-AA^\#=BB^\#-BXA^\#$ and $A(BB^\#-AA^\#)=YBBB^\#-A=A-A=0$. Therefore, $\R(BB^\#-AA^\#)\subseteq\R(B)\cap\N(A)$.
\end{proof}

\begin{proposition}\label{square_operator_matrix_form}
Let $A,B\in M_n^{\#}$ and $A\le^s B$. With respect to the
decomposition \eqref{dekomp_prostora},
\begin{equation}\label{reprezentacije_za_A_i_B_kada_je_Amanje_s_B}
A=\begin{bmatrix}A_1&0&0\\0&0&0\\0&0&0\end{bmatrix},\qquad
B=\begin{bmatrix}B_1&B_2&0\\B_3&B_4&0\\0&0&0\end{bmatrix},
\end{equation}
where $A_1:\R(A)\to \R(A)$ and
$M=\begin{bmatrix}B_1&B_2\\B_3&B_4\end{bmatrix}: \R(A)\oplus\left(\N(A)\cap\R(B)\right) \to \R(A)\oplus\left(\N(A)\cap\R(B)\right)$
are both invertible.
Moreover, $A\le^\square B$ if and only if $B_1=A_1$.
\end{proposition}
\begin{proof}
Assume that $A,B\in M_n^{\#}$ and $A\le^s B$. By Lemma~\ref{lema_space_order_decomposition}, we have decompositions \eqref{prva}, \eqref{druga} and \eqref{dekomp_prostora}. It follows that $A$ and $B$ have the representations \eqref{reprezentacije_za_A_i_B_kada_je_Amanje_s_B} with respect to \eqref{dekomp_prostora}, where $A_1=A|_{\R(A)}$ and $M=B|_{\R(B)}$ are invertible.
It is easy to check that $A^3=ABA$ if and only if $B_1=A_1$.
\end{proof}

We can rewrite Proposition~\ref{square_operator_matrix_form} in the following way.

\begin{corollary}\label{square_matrix_form}
Let $A,B\in M_n^\#$ and $a=\operatorname{rank}(A)$, $b=\operatorname{rank}(B)$. Then $A\le^\square B$ if and only if there exists an invertible matrix $K\in M_n$ such that
\begin{equation}\label{matrix_forms_eq}
A=K\begin{bmatrix} A_1&0&0\\0&0&0\\0&0&0\end{bmatrix} K^{-1} \quad \text{ and }\quad
B=K\begin{bmatrix} A_1&B_2&0\\B_3&B_4&0\\0&0&0\end{bmatrix} K^{-1},
\end{equation}
where $A_1\in M_a$ and $\begin{bmatrix} A_1&B_2\\B_3&B_4\end{bmatrix}\in M_b$ are both invertible.
\end{corollary}

\begin{remark}\label{remark_Schur_complement}
Denote $M:=\begin{bmatrix} A_1&B_2\\B_3&B_4\end{bmatrix}$ and let $S:=B_4-B_3A_1^{-1}B_2$ be the Schur complement of the block $A_1$ in $M$. Since $M\in M_b$ is invertible, $S$ is invertible, thus
\begin{equation}\label{Gutman_formula}
\operatorname{rank}(S)=b-a.
\end{equation}
Furthermore,
\begin{equation}\label{formula_za_M-1}
M^{-1}=\begin{bmatrix} A_1^{-1}+A_1^{-1}B_2 S^{-1} B_3 A_1^{-1} & -A_1^{-1}B_2 S^{-1} \\ -S^{-1}B_3 A_1^{-1} & S^{-1} \end{bmatrix}.
\end{equation}
It follows that
$$M^{-1}-\operatorname{diag}(A_1^{-1},0)
=\begin{bmatrix} A_1^{-1}B_2 \\ -I \end{bmatrix} S^{-1} \begin{bmatrix} B_3A_1^{-1} & -I \end{bmatrix}.$$
This is a product of three factors whose leftmost and rightmost factors have full column/row rank, respectively (due to the $-I$ block), and hence
\begin{equation}\label{Gutman_formula_posledica}
\operatorname{rank}\big(M^{-1}-\operatorname{diag}(A_1^{-1},0)\big)=\operatorname{rank}(S^{-1})=\operatorname{rank}(S).
\end{equation}
It is clear that
\begin{equation}\label{formule_za_grupne_inverze}
A^\#=K\operatorname{diag}(A_1^{-1},0,0)K^{-1} \text{ and } B^\#=K\operatorname{diag}(M^{-1},0)K^{-1}.
\end{equation}
It follows from \eqref{Gutman_formula}, \eqref{Gutman_formula_posledica} and \eqref{formule_za_grupne_inverze}
$$\operatorname{rank}(B^\#-A^\#)=b-a=\operatorname{rank}(B^\#)-\operatorname{rank}(A^\#)$$
which is equivalent to $A^\#\le^- B^\#$.
Thus, we obtain an alternative derivation of the result
$A\le^\square B \;\Rightarrow\; A^\#\le^- B^\#$ stated in Theorem~\ref{th_square_equiv_minus}.
\end{remark}

The following theorem establishes a connection between the reverse order law and the square and sharp partial orders.

\begin{theorem}\label{obrni_okreni}
Let $A,B\in M_n^\#$ and $a=\operatorname{rank}(A)$, $b=\operatorname{rank}(B)$. If $A\le^\square B$ then
$$(AB)^\#=B^\# A^\# \iff A\le^\# B.$$
\end{theorem}
\begin{proof}
Assume that $A\le^\square B$. Then, by Corollary~\ref{square_matrix_form}, \eqref{matrix_forms_eq} holds,
where $A_1$ and $M=\begin{bmatrix} A_1&B_2\\B_3&B_4\end{bmatrix}$ are both invertible.
Recall that $A\le^\# B$ if and only if $A^2=AB=BA$. Using this, it is easy to check that $A\le^\# B$ if and only if $B_2=0$ and $B_3=0$.
Let $M^{-1}=\begin{bmatrix}X_1&X_2\\X_3&X_4\end{bmatrix}$, where the explicit expressions of $X_i$ are given in \eqref{formula_za_M-1}.
Using \eqref{formule_za_grupne_inverze}, we obtain
\begin{align*}
(AB)(B^\#A^\#)&=A(BB^\#)A^\#=K\operatorname{diag}(I_a,0,0)K^{-1} \\
(B^\#A^\#)(AB)&=B^\#(A^\#A)B=K\begin{bmatrix} X_1A_1&X_1B_2&0\\X_3A_1&X_3B_2&0\\0&0&0 \end{bmatrix}K^{-1}.
\end{align*}
Assume that $(AB)^\#=B^\#A^\#$. Then $(AB)(B^\#A^\#)=(B^\#A^\#)(AB)$, hence $X_1A_1=I_a$, $X_1B_2=0$ and $X_3A_1=0$. It follows that $X_1=A_1^{-1}$, $B_2=0$ and $X_3=-S^{-1}B_3 A_1^{-1}=0$, so $B_3=0$. Therefore, $A\le^\# B$.
Conversely, if $A\le^\# B$ then $B_2=0$ and $B_3=0$, and hence, $B_4$ has to be invertible. Therefore, $AB=\operatorname{diag}(A_1^2,0,0)$ and $B^\#A^\#=\operatorname{diag}((A_1^{-1})^2,0,0)$, so $(AB)^\#=B^\# A^\#$.
\end{proof}

\subsection{Additive properties of the group inverse and related orders.}

The results obtained so far show that the square order has many properties analogous to those of the diamond order.
We now investigate the additive properties of the group inverse and show that its relationship with the sharp and square partial orders differs from that between the additive properties of the Moore–Penrose inverse and the diamond and star partial orders.
Hartwig and Styan showed in \cite[Theorem 2]{Hartwig_Styan} that
$$A\le^- B \text{ and } (B-A)^\dag=B^\dag-A^\dag \iff A\le^* B.$$
A purely algebraic proof of this equivalence was given by Baksalary \cite{Baksalary_aditivity_MP}, and a geometric proof in the setting of bounded linear operators on Hilbert spaces was given by Fongi and Gonzalez \cite{Fongi_Gonzalez}.
As a consequence, if $(B-A)^\dag=B^\dag-A^\dag$ then
$$A\le^* B \iff A\le^\diamondsuit B \iff A\le^-B.$$

\begin{remark}
Let $A,B\in M_n^{\#}$. In general,
$$A\le^- B  \text{ and } (B-A)^\#=B^\#-A^\# \nimplies  A\le^{\#}_{\ell}B.$$
Indeed, let $A=\operatorname{diag}(1,0)$ and
$B=\begin{bmatrix}1&1\\0&-1\end{bmatrix}$. Then
$\operatorname{rank}(B)=2=\operatorname{rank}(A)+\operatorname{rank}(B-A)$,
so $A\le^{-}B$. Since $B^2=I_2$ and $(B-A)^3=B-A$, we have $B^{\#}=B$, $A^{\#}=A$ and $(B-A)^{\#}=B-A=B^{\#}-A^{\#}$.
However, $AB=\begin{bmatrix}1&1\\0&0\end{bmatrix}\neq A^2$. Replacing $B$ by $B^{T}$, we see that $A\le^\#_rB$ is not implied either.
\end{remark}

\begin{theorem}\label{additivity_group_th}
Let $A,B\in M_n^\#$. Then
$$A\le^- B \text{ and } (B-A)^\#=B^\#-A^\#$$
if and only if
$$A^3=ABA=AB^2=B^2A \text{ and } (B-A)^3=B^3-A^3.$$
\end{theorem}
\begin{proof}
Assume that $A\le^- B$ and $(B-A)^\#=B^\#-A^\#$. Then $\r(B-A)=\r(B^\#-A^\#)$. By Corollary~\ref{square_preko_rangova}(iv),
we conclude that $A\le^\square B$. Thus, $A^3=ABA$ and $A\le^s B$. Therefore, $A=ABA^\#=A^\#BA=AB^\#B=BB^\#A$. From $(B-A)^\#=B^\#-A^\#$, we obtain
$(B-A)^2(B^\#-A^\#)=B-A$
which reduces to
$$B^2A^\#+BAB^\#=BAA^\#+A^2B^\#.$$
Multiplying the last equation by $BA$ from the right side, we obtain $A^3=B^2A$. We can obtain in a similar way that $A^3=AB^2$.
From $A\le^- B$ and $(B-A)^\#=B^\#-A^\#$, it follows that $B-A\le^- B$ and $(B-(B-A))^\#=B^\#-(B-A)^\#$. By the result established above, it follows that $(B-A)^3=B^2(B-A)$. Thus, $(B-A)^3=B^3-A^3$.

Assume now that $A^3=ABA=AB^2=B^2A$ and $(B-A)^3=B^3-A^3$ and let $a=\r(A)$, $b=\r(B)$. Then
$\R(A)=\R(A^3)=\R(B^2A)\subseteq\R(B)$ and $\N(B)\subseteq\N(AB^2)=\N(A^3)=\N(A)$, so
$A\le^s B$. Therefore, $A\le^\square B$. Let us first show that $A\le^-B$. If $G\in B\{1\}$ then
\begin{align*}
AGA&=(A^\#)^2A^3GA^3(A^\#)^2=(A^\#)^2AB^2GB^2A(A^\#)^2\\
&=(A^\#)^2AB^3A(A^\#)^2=(A^\#)^2A^3BA(A^\#)^2=(A^\#)^2A^5(A^\#)^2=A
\end{align*}
Therefore, $B\{1\}\subseteq A\{1\}$, so $A\le^- B$. By Corollary~\ref{square_matrix_form},
$$A=K\begin{bmatrix} A_1&0&0\\0&0&0\\0&0&0\end{bmatrix} K^{-1} \quad \text{ and }\quad
B=K\begin{bmatrix} A_1&B_2&0\\B_3&B_4&0\\0&0&0\end{bmatrix} K^{-1},$$
for some invertible matrices $K\in M_n$, $A_1\in M_a$ and $M=\begin{bmatrix} A_1&B_2\\B_3&B_4\end{bmatrix}\in M_b$.
We have
$$M^{-1}=\begin{bmatrix} X_1&X_2\\X_3&X_4\end{bmatrix}=
\begin{bmatrix} A_1^{-1}+A_1^{-1}B_2 S^{-1} B_3 A_1^{-1} & -A_1^{-1}B_2 S^{-1} \\ -S^{-1}B_3 A_1^{-1} & S^{-1} \end{bmatrix}$$
where $S=B_4-B_3A_1^{-1}B_2$. The condition $A^3=AB^2$ yields $AB^\#=A^\#B$. By \eqref{formule_za_grupne_inverze}, it reduces to $A_1X_1=I_a$ and $A_1X_2=A_1^{-1}B_2$. It follows that,
$X_1=A_1^{-1}$, so $B_2S^{-1}B_3=0$, and $-B_2S^{-1}=A_1^{-1}B_2$. Multiplying the last equation by $B_3$ from the right, we obtain $B_2B_3=0$.
Therefore,
$$B_2=-A_1^{-1}B_2S=-A_1^{-1}B_2B_4.$$
Similarly, from $A^3=B^2A$, we conclude that
$$-S^{-1}B_3=B_3A_1^{-1}, \quad B_3=-B_4B_3A_1^{-1}.$$
From $A^3=ABA=AB^2=B^2A$, the condition $(B-A)^3=B^3-A^3$ yields $A^3+BAB=A^2B+BA^2$ which reduces to
$B_3A_1B_2=0.$
As $B_3A_1^nB_2=-B_4B_3A_1^{n-1}B_2$, we conclude by induction on $n$ that
\begin{equation}\label{CHT}
B_3A_1^nB_2=0, \; n\geq 1.
\end{equation}
Let $p(x)=x^a+c_{a-1}x^{a-1}+\cdots +c_1x+c_0$ be the characteristic polynomial of $A_1$. Since $A_1$ is invertible, $c_0\neq 0$. From the Cayley--Hamilton theorem and \eqref{CHT}, we obtain that $0=B_3p(A_1)B_2=c_0B_3B_2$, that is,
$$B_3A_1^nB_2=0, \; n\geq 0.$$
Since $A_1^{-1}$ is a polynomial in $A_1$, it follows that
$B_3A_1^{-1}B_2=0$, so $S=B_4$, and thus $B_2B_4^{-1}B_3=0$. Therefore,
$$B^\#=K\begin{bmatrix} A_1^{-1}&B_2B_4^{-2}&0\\B_4^{-2}B_3&B_4^{-1}&0\\0&0&0\end{bmatrix}K^{-1}.$$
Since $B_2B_4^{-1}=-A_1^{-1}B_2$, we have $B_2B_4^{-n}B_3=(-1)^nA_1^{-n}B_2B_3=0$ for $n\ge0$. Using this and $B_3B_2=0$, it is easy to check that
$$(B-A)(B^\#-A^\#)=(B^\#-A^\#)(B-A)=K\begin{bmatrix}0&B_2B_4^{-1}&0\\ B_4^{-1}B_3&I_{b-a}&0\\ 0&0&0\end{bmatrix}K^{-1},$$
$(B-A)(B^\#-A^\#)(B-A)=B-A$ and $(B^\#-A^\#)(B-A)(B^\#-A^\#)=B^\#-A^\#$.
It follows that $B-A\in M_n^\#$ and $(B-A)^\#=B^\#-A^\#$.
\end{proof}

\begin{remark}\label{remark_pomocna_pu}
Let $A,B\in M_n^{\#}$. The conditions used in Theorem~\ref{additivity_group_th} lead to two pre-orders on $M_n^{\#}$:
\begin{align*}
A\le^s B \text{ and } A^3=AB^2 & \iff A^2\le^{\#}_{\ell}B^2, \\
A\le^s B \text{ and } A^3=B^2A &\iff A^2\le^{\#}_{r}B^2 .
\end{align*}
Both relations are reflexive and transitive, but not antisymmetric: for $B=-A$ we have $A^2=B^2$, and hence each of them relates $A$ to $B$ and $B$ to $A$.
\end{remark}

We can define the following relations on $M_n^\#$.
\begin{definition}
Let $A,B\in M_n^\#$. Write
\begin{align*}
A\le^\boxplus_{\ell} &B  \iff A\le^s B \text{ and } A^3=ABA=AB^2 \\
A\le^{\boxplus}_r &B  \iff A\le^s B \text{ and } A^3=ABA=B^2A \\
A\le^\boxplus &B  \iff A^3=ABA=AB^2=B^2A.
\end{align*}
\end{definition}
In the definition of $\le^{\boxplus}$ the condition $A\le^s B$ is omitted, since it follows from $A^3=AB^2=B^2A$.

\begin{proposition}\label{prop_boxplus}
The relations $\le^{\boxplus}_{\ell}$, $\le^{\boxplus}_{r}$ and $\le^{\boxplus}$ are partial orders on $M_n^{\#}$. Moreover, each of them implies both $\le^{\square}$ and $\le^{-}$.
\end{proposition}
\begin{proof}
Since $\le^{\boxplus}_{\ell}$, $\le^{\boxplus}_{r}$ and $\le^{\boxplus}$ are intersections of the partial order $\le^{\square}$ (Corollary~\ref{cor_square_je_pu}) with the pre-orders from Remark~\ref{remark_pomocna_pu}, they are partial orders on $M_n^{\#}$.
Assume that $A\le^\boxplus_{\ell} B$. Then $A=BX=YB$ for some $X,Y\in M_n$, and for arbitrary $B^-\in B\{1\}$, we have
\begin{align*}
AB^-A&=YBB^-BX=YBX=AX=(A^\#)^2A^3X\\
&=(A^\#)^2AB^2X=(A^\#)^2ABA=(A^\#)^2A^3=A.
\end{align*}
Therefore, $B\{1\}\subseteq A\{1\}$, and thus, $A\le^-B$. Similarly, $A\le^\boxplus_r B\, \Rightarrow \, A\le^-B$, and so
$A\le^\boxplus B\, \Rightarrow \, A\le^-B$.
\end{proof}

We can rewrite Theorem~\ref{additivity_group_th} in the following way.
\begin{corollary}\label{additivity_group_cor}
Let $A,B\in M_n^\#$. Then
$$A\le^- B \text{ and } (B-A)^\#=B^\#-A^\# \iff A\le^\boxplus B \text{ and } B-A\le^\boxplus B.$$
\end{corollary}
\begin{proof}
$\Rightarrow$: By Theorem~\ref{additivity_group_th}, $A\le^\boxplus B$. Moreover, $B-A\in M_n^\#$, $B-A\le^- B$ and $(B-(B-A))^\#=B^\#-(B-A)^\#$, so Theorem~\ref{additivity_group_th} applied to $B-A$ and $B$ gives $B-A\le^\boxplus B$.

$\Leftarrow$: From $B-A\le^\boxplus B$ and $B^2A=A^3$ we obtain $(B-A)^3=B^2(B-A)=B^3-A^3$, so the conditions of
Theorem~\ref{additivity_group_th} are satisfied.
\end{proof}

\subsection{Comparison with other orders.}

We now compare $\le^\square$, $\le^\boxplus_\ell$, $\le^\boxplus_r$ and
$\le^\boxplus$ with the minus order and the sharp-type orders.

\begin{proposition}\label{prop:comparison}
Let $A,B\in M_n^\#$. Then
$A\le^\#_\ell B\Rightarrow A\le^\boxplus_\ell B$,
$A\le^\#_r B\Rightarrow A\le^\boxplus_r B$ and
$A\le^\# B\Rightarrow A\le^\boxplus B$.
\end{proposition}

\begin{proof}
If $A^2=AB$, then $ABA=A^3$ and $AB^2=A^2B=A^3$; the case $A^2=BA$ is
symmetric, and $\le^\#$ is the intersection of $\le^\#_\ell$ and $\le^\#_r$.
\end{proof}

\begin{figure}[htbp]
\centering
\begin{tikzpicture}[
    x=1.55cm, y=1.35cm,
    ord/.style={draw, rounded corners=2pt, minimum width=1.25cm,
                minimum height=0.62cm, inner sep=2pt},
    new/.style={ord, fill=black!18},
    imp/.style={-{Stealth[length=2.2mm]}, semithick},
  ]
  \node[ord] (sh)  at ( 0, 4) {$\le^{\#}$};
  \node[ord] (shl) at (-2, 3) {$\le^{\#}_{\ell}$};
  \node[new] (bp)  at ( 0, 3) {$\le^{\boxplus}$};
  \node[ord] (shr) at ( 2, 3) {$\le^{\#}_{r}$};
  \node[new] (bpl) at (-1, 2) {$\le^{\boxplus}_{\ell}$};
  \node[new] (bpr) at ( 1, 2) {$\le^{\boxplus}_{r}$};
  \node[new] (sq)  at (-1, 1) {$\le^{\square}$};
  \node[ord] (mi)  at ( 1, 1) {$\le^{-}$};

  \draw[imp] (sh)  -- (shl);
  \draw[imp] (sh)  -- (bp);
  \draw[imp] (sh)  -- (shr);
  \draw[imp] (shl) -- (bpl);
  \draw[imp] (bp)  -- (bpl);
  \draw[imp] (bp)  -- (bpr);
  \draw[imp] (shr) -- (bpr);
  \draw[imp] (bpl) -- (sq);
  \draw[imp] (bpl) -- (mi);
  \draw[imp] (bpr) -- (sq);
  \draw[imp] (bpr) -- (mi);
\end{tikzpicture}
\caption{Implications between the partial orders on $M_n^\#$ studied in
Section~2; the new orders are shaded.}
\label{dijagram1}
\end{figure}
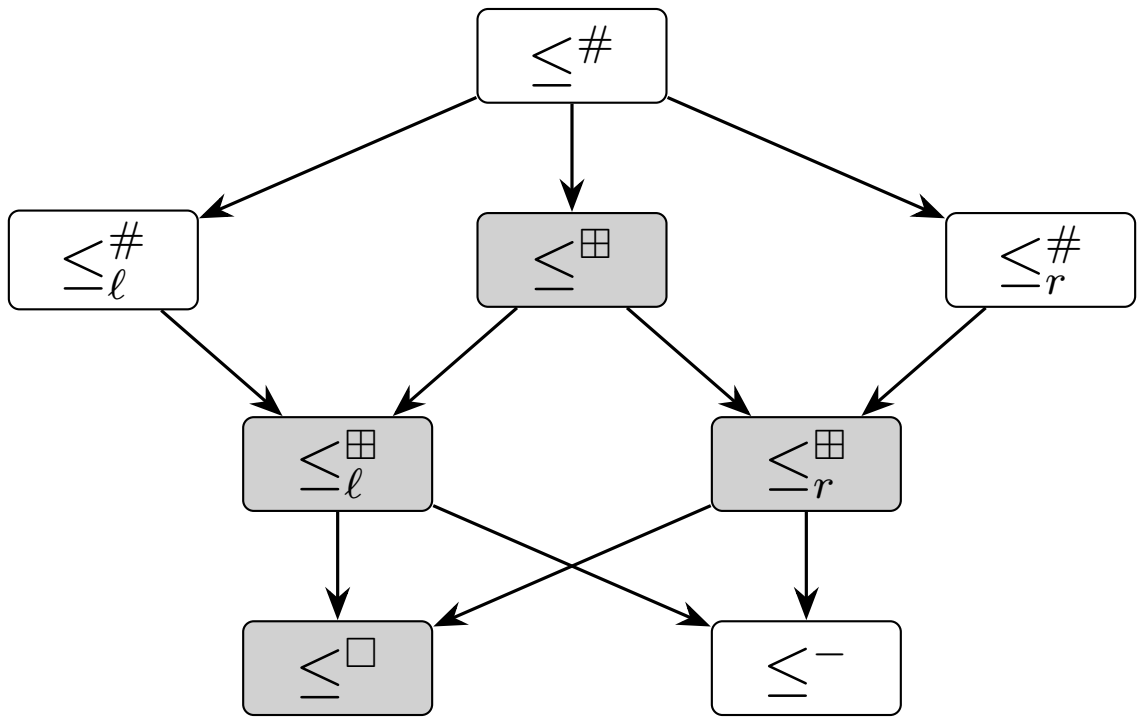

Together with Proposition~\ref{prop_boxplus}, this gives all implications
shown in Figure~\ref{dijagram1}. In particular, each of the one-sided sharp
orders implies both $\le^\square$ and $\le^-$.
Remark~\ref{remark_counter} shows that no other implications hold.

\begin{remark}\label{remark_counter}
Let $A,B\in M_n^\#$. In general,
\begin{enumerate}[label=(\roman*)]
\item $A\le^\square B \nimplies A\le^- B$ and
      $A\le^\square B \nimplies A^\#\le^\square B^\#$;
\item $A\le^- B \nimplies A\le^\square B$;
\item $A\le^\#_r B \nimplies A\le^\boxplus_\ell B$ and
      $A\le^\#_\ell B \nimplies A\le^\boxplus_r B$;
\item $A\le^\boxplus B \nimplies A\le^\#_\ell B$ and
      $A\le^\boxplus B \nimplies A\le^\#_r B$;
\item $A\le^\square B \nimplies B-A\in M_n^\#$.
\end{enumerate}
Consequently, no implications between the orders in
Figure~\ref{dijagram1} hold other than those obtained by composing the
arrows.
\end{remark}

\begin{proof}
In (i)--(iii), let $A=\operatorname{diag}(1,0)$.

(i) Let $B=\begin{bmatrix}1&1\\1&0\end{bmatrix}$. Then $A\le^\square B$, since $B$ is invertible and $ABA=A=A^3$. However,
$\operatorname{rank}(B-A)=2$, so $A\le^- B$ does not hold. By Theorem~\ref{th_square_equiv_minus}, $A^\#\le^\square B^\#$ if and only if $A\le^- B$, so $A^\#\le^\square B^\#$ does not hold either.

(ii) Let $B=\begin{bmatrix}2&1\\1&1\end{bmatrix}$. Then $\operatorname{rank}(B-A)=1$, so $A\le^- B$, but $ABA=2A\neq A^3$.

(iii) Let $B=\begin{bmatrix}1&1\\0&1\end{bmatrix}$. Then $A^2=BA$, so $A\le^\#_r B$, but $AB^2\neq A^3$, so $A\le^\boxplus_\ell B$ does not hold. For $B^T$ the roles of $\ell$ and $r$ are interchanged.

(iv) Let $A=\operatorname{diag}(1,0,0)$ and $B=\begin{bmatrix}1&0&1\\1&-1&1\\0&0&-1\end{bmatrix}$. Then $B$ is invertible and $A^3=ABA=AB^2=B^2A$, so $A\le^\boxplus B$, but $AB\ne A^2$ and $BA\ne A^2$.

(v) Let
$A=\begin{bmatrix}1&0&0\\0&1&0\\0&0&0\end{bmatrix}$ and
$B=\begin{bmatrix}1&0&1\\0&1&0\\0&1&1\end{bmatrix}$. Then $A\le^\square B$, since $B$ is invertible and $ABA=A=A^3$, but
$\operatorname{rank}(B-A)=2\neq1=\operatorname{rank}((B-A)^2)$.
\end{proof}

\section{Unification of partial orders arising from the plus partial order}

\subsection{The plus order and g-maps.}

Motivated by the minus, star, and diamond partial orders, Arias and Maestripieri recently introduced the \emph{plus partial order}, denoted by $\le^+$, in \cite{Arias_Maestripieri}, which is closely related to the relation $\le^\square$.
\begin{definition}[\cite{Arias_Maestripieri}]\label{plus_order_def}
If $A,B\in M_{m\times n}$ then
$A\le^+ B$ if and only if $A\le^s B$ and there exist idempotents $P\in E(M_m)$ and $Q\in E(M_n)$ such that $A=PBQ$.
In such a case, without loss of generality we can assume that $\R(P)=\R(A)$ and $\R(Q^*)=\R(A^*)$.
\end{definition}
The relation $\le^+$ is a partial order on $M_{m\times n}$.
It was originally defined in \cite{Arias_Maestripieri} for bounded operators on Hilbert spaces and it was extended to
Rickart $*$-rings in \cite{MarovtMosic_plus}.

\begin{proposition}[Proposition 4.8 in \cite{Arias_Maestripieri}]\label{plus_<_sledi_plus_preceq}
Let $A,B\in M_{m\times n}$. If $A\le^+ B$ then there exist $A^-\in A\{1,2\}$ and $B^-\in B\{1,2\}$ such that $A^-\le^- B^-$.
\end{proposition}

From Proposition~\ref{square_sledi_plus}, it follows that
$$A\le^\square B \implies A\le^+ B.$$
Recall the following well-known characterization of the diamond partial order:
\begin{equation}
\begin{aligned}\label{diamond_kao_plus}
&A\le^\diamondsuit B \iff \\
&A\le^s B \text{ and } A=PBQ, \text{ for some } P\in E^*(M_m), \; Q\in E^*(M_n),
\end{aligned}
\end{equation}
where without loss of generality we can choose $P=AA^\dag$ and $Q=A^\dag A$.
Therefore,
$$A\le^\diamondsuit B \implies A\le^+ B.$$

Motivated by the g-based partial orders (i.e., partial orders based on generalized inverses) introduced by Mitra in \cite{Unified_theory}, we go further and develop a unified framework for partial orders arising from the plus partial order. In this framework, the plus order plays a role analogous to that of the minus order in the theory of g-based orders.

We briefly introduce Mitra's terminology given in \cite{Unified_theory} and in Chapter 7 of \cite{knjiga}.
Let $\mathcal{P}(M_{n\times m})$ denote the class of all subsets of $M_{n\times m}$.
\begin{definition}
A \emph{g-map} is a map
$\mathcal{G}:M_{m\times n}\to\mathcal{P}(M_{n\times m})$
such that for each $A\in M_{m\times n}$, $\mathcal{G}(A)$ is a certain subset (possibly empty) of $A\{1\}$.
The set
$$\Omega_{\mathcal{G}}=\{A\in M_{m\times n} : \mathcal{G}(A)\neq\emptyset\}$$
is called the \emph{support} of $\mathcal{G}$. We write
$$\G_r(A)=\G(A)\cap A\{2\}.$$
\end{definition}
\begin{definition}
We say that $\G$ is \emph{semi-complete} if for each $A\in \Omega_\G$, $\{A^-AA^- : A^-\in \G(A)\}\subseteq \G(A)$.
\end{definition}
Mitra's concept of unification of g-based partial orders is as follows.

\begin{definition}\label{Mitra_g_order}
Let $\G:M_{m\times n}\to\mathcal{P}(M_{n\times m})$ be a g-map and $A,B\in M_{m\times n}$. Then $A\le^{\mathcal{G}} B$ if there exists $A^-\in \G(A)$ such that $AA^-=BA^-$ and $A^-A=A^-B$.
\end{definition}

The following lemma is related to the fact that a $\{1,2\}$-inverse of $A$
is uniquely determined by its range and null space
(see \cite[Chapter~2, Section~6]{BIG}).

\begin{lemma}\label{remark_P,Q}
Let $A\in M_{m\times n}$, $P\in E(M_m)$ and $Q\in E(M_n)$ be such that $\R(P)=\R(A)$ and $\N(Q)=\N(A)$, and let $X\in A\{1\}$ be arbitrary. Then $G=QXP$ is the unique matrix in $A\{1,2\}$ such that $AG=P$ and $GA=Q$. Consequently,
\begin{equation}
\begin{aligned}\label{P_i_Q_preko_A-}
&\{(P,Q): P\in E(M_m),\, Q\in E(M_n), \, \R(P)=\R(A), \, \N(Q)=\N(A) \} \\
&=\{(AA^-,A^-A) : A^-\in A\{1\}\}.
\end{aligned}
\end{equation}
\end{lemma}

\begin{proof}
Since $\R(P)=\R(A)$ and $P$ is idempotent, we have $PA=A$ and $P=AY$ for some $Y$; similarly, $AQ=A$ and $Q=ZA$ for some $Z$. Hence $AG=AQXP=AXAY=AY=P$, $GA=QXPA=ZAXA=ZA=Q$, $AGA=PA=A$ and $GAG=QG=G$, so $G\in A\{1,2\}$. If $G'\in A\{1,2\}$ also satisfies $AG'=P$ and $G'A=Q$, then
$$G'=G'AG'=G'AGAG'=QGP=(GA)G(AG)=G.$$
Finally, for $A^-\in A\{1\}$ the matrices $AA^-$ and $A^-A$ are idempotents with $\R(AA^-)=\R(A)$ and $\N(A^-A)=\N(A)$, which proves ``$\supseteq$'' in \eqref{P_i_Q_preko_A-}, and ``$\subseteq$'' follows from $(P,Q)=(AG,GA)$.
\end{proof}

\subsection{The relations $\le^{\G}_{+}$ and $\preceq^{\G}_{+}$.}

We are now in a position to introduce two new classes of relations.

\begin{definition}\label{def<}
Let $\G:M_{m\times n}\to\mathcal{P}(M_{n\times m})$ be a g-map and $A,B\in M_{m\times n}$. We say that $A\le^{\G}_{+} B$ if
$A\le^s B$ and there exists
$A^-\in \G(A)$ such that $A=AA^-BA^-A$.
\end{definition}

If $A^-\in A\{1\}$ then $A^-AA^-\in A\{1,2\}$. Thus, if $\G$ is a semi-complete g-map, then
$$A\le^{\G}_{+} B \iff A\le^s B \text{ and } \exists A^-\in \G_r(A) \text{ such that } A^-=A^-BA^-.$$

\begin{definition}\label{def_preceq}
Let $\G:M_{m\times n}\to\mathcal{P}(M_{n\times m})$ be a g-map and $A,B\in M_{m\times n}$. We say that $A\preceq^{\mathcal{G}}_+ B$ if $A\le^s B$ and there exist
$A^-\in \mathcal{G}_r(A)$ and $B^-\in \mathcal{G}_r(B)$ such that $A^-\le^- B^-$.
\end{definition}

If $B^-\in \mathcal{G}_r(B)$ then $B^-\in B\{1,2\}$, so $B\in B^-\{1\}$. Therefore, by Theorem~\ref{minus_karakterizacije}, it follows that
\begin{equation}\label{preceq_karakterizacija}
A\preceq^{\G}_{+} B \iff A\le^s B \text{ and } A^-=A^-BA^-=A^-BB^-=B^-BA^-,
\end{equation}
for some $A^-\in \G_r(A)$ and $B^-\in \mathcal{G}_r(B)$.

We call $\le^{\G}_{+}$ and $\preceq^{\G}_{+}$ the \emph{g-plus order relation} and the \emph{strong g-plus order relation} induced by $\G$, respectively.

\begin{remark}\label{RemarkGr}
Note that $\preceq^{\G}_{+}$ depends only on $\G_r$. Moreover, if $\G$ is semi-complete, then $\le^{\G}_{+}$ also depends only on $\G_r$. Indeed, if $A\le^s B$ and $A=AA^-BA^-A$ for some $A^-\in\G(A)$, then $H:=A^-AA^-\in\G_r(A)$ and $AHBHA=A$. Thus, for a semi-complete g-map $\G$, the relations $\le^{\G}_{+}$ and $\preceq^{\G}_{+}$ coincide with $\le^{\G_r}_+$ and $\preceq^{\G_r}_+$, respectively, and $\Omega_\G=\Omega_{\G_r}$.
\end{remark}

\begin{proposition}\label{preceq_sledi_<}
Let $\G:M_{m\times n}\to\mathcal{P}(M_{n\times m})$ be a g-map and $A,B\in M_{m\times n}$. If $A\preceq^{\G}_{+} B$ then $A\le^{\G}_{+} B$.
\end{proposition}
\begin{proof}
Assume that $A\preceq^{\G}_{+} B$. Then $A\le^s B$ and $A^-\le^- B^-$ for some $A^-\in \G_r(A)$ and $B^-\in \G_r(B)$.
From \eqref{preceq_karakterizacija}, we obtain $A=AA^-A=AA^-BA^-A$, so $A\le^{\G}_{+} B$.
\end{proof}

Thus, the strong g-plus order relation implies the g-plus order relation, which justifies the term \emph{strong}.

\begin{example}\label{primer_G}
Relations $\le^{\G}_{+}$ and $\preceq^{\G}_{+}$ defined in Definitions~\ref{def<} and \ref{def_preceq} include the plus, diamond, and square partial orders as special cases.
\begin{enumerate}[label=(\roman*)]
\item Let $\G(A)=A\{1\}$. Then $\G$ is semi-complete and $\G_r(A)=A\{1,2\}$. Hence, by Remark~\ref{RemarkGr}, the relations $\le^{\G}_{+}$ and $\preceq^{\G}_{+}$ do not change if $A\{1\}$ is replaced by $A\{1,2\}$. By \eqref{P_i_Q_preko_A-}, $\le^{\G}_{+}$ is the plus partial order $\le^+$ from Definition~\ref{plus_order_def}.
    Moreover, by Proposition~\ref{plus_<_sledi_plus_preceq}, $A\le^+B$ implies $A\preceq^{\G}_{+}B$, and by Proposition~\ref{preceq_sledi_<}, $A\preceq^{\G}_{+}B$ implies $A\le^{\G}_{+}B$. That is, if $\G(A)=A\{1\}$ or $\G(A)=A\{1,2\}$, then
    $$A\le^{\G}_{+} B \iff A\preceq^{\G}_{+} B \iff A\le^+ B .$$
\item Let $\G$ be defined by $\G(A)=\{A^\dag\}$. From \eqref{diamond_kao_plus} and \eqref{diamond_kao_preceq}, it follows that for $A,B\in M_{m\times n}$,
    $$A\le^{\G}_{+} B \iff A\preceq^{\G}_{+} B \iff A\le^\diamondsuit B.$$
\item Let a g-map $\G:M_{n}\to\mathcal{P}(M_{n})$ be defined by $\G(A)=\{A^\#\}$, $\Omega_{\mathcal{G}}=M_n^\#$. From Proposition~\ref{square_sledi_plus} and Theorem~\ref{th_square_equiv_minus}, it follows that for $A,B\in M_n^\#$,
    $$A\le^{\G}_{+} B \iff A\preceq^{\G}_{+} B \iff A\le^\square B.$$
\end{enumerate}
These are the cases $S=\emptyset$, $\{2,3,4\}$ and $\{2,5\}$ of Theorem~\ref{teorema_klasifikacija}.
\end{example}

For $A\in M_{m\times n}$ we write
\begin{align*}
&\PP(A)=\{P\in E(M_m):\mathcal R(P)=\mathcal R(A)\},\\
&\QQ(A)=\{Q\in E(M_n):\mathcal N(Q)=\mathcal N(A)\},
\end{align*}
and, for a g-map $\G$ and $A\in\Omega_\G$,
\[
\Pi_\G(A)=\{(AA^-,A^-A): A^-\in\G(A)\}\subseteq \PP(A)\times\QQ(A).
\]
By Lemma~\ref{remark_P,Q}, for every $(P,Q)\in\PP(A)\times\QQ(A)$ there exists a unique $A^-\in A\{1,2\}$ such that $AA^-=P$ and $A^-A=Q$, namely $A^-=QXP$, where $X\in A\{1\}$ is arbitrary.

\begin{proposition}\label{proposition_pairs}
Let $\G:M_{m\times n}\to\mathcal{P}(M_{n\times m})$ be a g-map, $A\in\Omega_\G$ and $B\in M_{m\times n}$.
\begin{enumerate}[label=(\roman*)]
\item $A\le^{\G}_{+}B$ if and only if $A\le^s B$ and $A=PBQ$ for some
      $(P,Q)\in\Pi_\G(A)$.
\item If $\G$ is semi-complete, then $A\le^{\G}B$ if and only if
      $A=PB=BQ$ for some $(P,Q)\in\Pi_\G(A)$.
\end{enumerate}
In particular, both relations depend only on $\Pi_\G$.
\end{proposition}

\begin{proof}
(i) follows from $AA^-BA^-A=(AA^-)B(A^-A)$.

(ii) Assume that $AA^-=BA^-$ and $A^-A=A^-B$ for some $A^-\in\G(A)$.
Let $P=AA^-$ and $Q=A^-A$. Then $(P,Q)\in \Pi_\G(A)$ and $A=PB=BQ$.
Conversely, assume that $A=PB=BQ$ for some $(P,Q)\in\Pi_\G(A)$, that is, $A=AA^-B=BA^-A$ for some $A^-\in\G(A)$. Since $\G$ is semi-complete, $H=A^-AA^-\in\G(A)$ and $AH=BH$, $HA=HB$, i.e. $A\le^\G B$.
\end{proof}

Part (ii) is closely related to \cite[Theorem 7.5.1]{knjiga}, where the projectors are described in terms of $\R(G)$ and $\N(G)$, $G\in\G(A)$.

Since $AA^-$ and $A^-A$ are idempotents, Definition~\ref{plus_order_def} shows that $A\le^{\G}_{+}B$ implies $A\le^{+}B$. Together with Proposition~\ref{preceq_sledi_<}, this gives
$$A\preceq^{\G}_{+}B \implies A\le^{\G}_{+}B \implies A\le^{+}B .$$
In particular, the plus order, obtained for $\G(A)=A\{1\}$ (Example~\ref{primer_G}(i)), is the largest of the relations $\le^{\G}_{+}$.

\begin{proposition}\label{reflanti}
Let $\G:M_{m\times n}\to\mathcal{P}(M_{n\times m})$ be a g-map. The relation $\le^{\G}_{+}$ is reflexive and antisymmetric on $\Omega_\G$. If $\G$ is semi-complete, then the relation $\preceq^{\G}_{+}$ is reflexive and antisymmetric on $\Omega_\G$.
\end{proposition}
\begin{proof}
The reflexivity of $\le^{\G}_{+}$ is obvious. Assume that $A\le^{\G}_{+} B$ and $B\le^{\G}_{+} A$.
Since $A\le^{\G}_{+} B$, there exists $A^-\in \G(A)$ such that $A=AA^-BA^-A$.
From $B\le^s A$, we have $B=AA^-B=BA^-A$, so $A=B$.
Therefore, $\le^{\G}_{+}$ is antisymmetric.
Assume now that $\G$ is semi-complete and let $A\in\Omega_\G$. For $A^-\in \G(A)$, we know that $G:=A^-AA^-\in \G_r(A)$. From $G\le^-G$, it follows that $A\preceq^{\G}_{+} A$. Thus, $\preceq^{\G}_{+}$ is reflexive.
It follows from  Proposition~\ref{preceq_sledi_<} that $\preceq^{\G}_{+}$ is also antisymmetric on $\Omega_{\G}$.
\end{proof}

\begin{remark}
We could have defined the relation $\preceq^{\G}_{+}$ in the following way: $A\preceq^{\G}_{+} B$ if and only if $A\le^s B$ and there exist $A^-\in \G(A)$ and $B^-\in \G(B)$ such that $A^- \le^- B^-$. However, the relation defined in this way would not be antisymmetric. The following example demonstrates this. Let
$A=\begin{bmatrix}1&0\\0&0\end{bmatrix}$, $B=\begin{bmatrix}2&0\\0&0\end{bmatrix}$,
$G_1=A$, $G_2=\begin{bmatrix}1&1\\1&2\end{bmatrix}$, $H_1=\begin{bmatrix}0.5&1\\1&-2\end{bmatrix}$ and $H_2=\begin{bmatrix}0.5&0\\0&0\end{bmatrix}$.
It is readily verified that $G_1,G_2\in A\{1\}$ and $H_1,H_2\in B\{1\}$, but $G_2\notin A\{2\}$ and $H_1\notin B\{2\}$.
Let $\G(A)=\{G_1,G_2\}$ and $\G(B)=\{H_1,H_2\}$.
Then $B=2A$, so $A\le^s B$ and $B\le^s A$.
Moreover, $\operatorname{rank}(H_1-G_1)=1=2-1=\operatorname{rank}(H_1)-\operatorname{rank}(G_1)$ and $\operatorname{rank}(G_2-H_2)=1=2-1=\operatorname{rank}(G_2)-\operatorname{rank}(H_2)$. Therefore, $G_1\le^-H_1$ and $H_2\le^-G_2$, which implies that $A\preceq^{\G}_{+} B$ and $B\preceq^{\G}_{+} A$. However, $A\neq B$, and hence $\preceq^{\G}_{+}$ is not antisymmetric.

On the other hand, $AG_1BG_1A=B\neq A$ and $AG_2BG_2A=B\neq A$, so that $A\le^{\G}_{+} B$ does not hold. Consequently, Proposition~\ref{preceq_sledi_<} would no longer hold.

These are some of the reasons for requiring the conditions $A^-\in \G_r(A)$ and $B^-\in \G_r(B)$ in Definition~\ref{def_preceq}.
\end{remark}

\begin{proposition}\label{<_ne_sledi_preceq}
Let $\G:M_{m\times n}\to\mathcal{P}(M_{n\times m})$ be a g-map.
\begin{enumerate}[label=(\roman*)]
\item If $A\in\Omega_{\G}$ and $B\in M_{m\times n}$, then $A\le^{\G} B$ implies $A\le^{\G}_{+} B$.
\item In general, $A\le^{\G}B$ does not imply $A\preceq^{\G}_{+}B$; that is,
      there exist a g-map $\G$ and $A,B\in\Omega_{\G}$ such that
      $A\le^{\G}B$ (and hence $A\le^{\G}_{+}B$) but $A\preceq^{\G}_{+}B$ does not hold.
\end{enumerate}
\end{proposition}
\begin{proof}
(i)
Let $A\le^{\G} B$. Then $AA^-=BA^-$ and $A^-A=A^-B$ for some $A^-\in \G(A)$. We obtain $A=AA^-A=BA^-A=AA^-B$. Hence, $A\le^s B$ and $A=AA^-BA^-A$. Thus, $A\le^{\G}_{+} B$.

(ii)
Let
$A=\begin{bmatrix} 1&0&0\\0&0&0\\0&0&0\end{bmatrix}$, $B=\begin{bmatrix} 1&0&0\\0&1&0\\0&0&0\end{bmatrix}$, $A^-=A$ and
$B^-=\begin{bmatrix} 1&0&1\\0&1&1\\0&0&0\end{bmatrix}$. Then $A\le^s B$ and $A^-\in A\{1,2\}$, $B^-\in B\{1,2\}$. Let $\G$ be a g-map such that $\G(A)=\{A^-\}$ and $\G(B)=\{B^-\}$; then $\G(A)=\G_r(A)$ and $\G(B)=\G_r(B)$. We have $AA^-=BA^-$ and $A^-A=A^-B$, so $A\le^{\G} B$, and by (i) $A\le^{\G}_{+} B$.
But $A^-\neq A^-BB^-$, so, by \eqref{preceq_karakterizacija} (the only possible choices being $A^-$ and $B^-$), $A^-\le^- B^-$ does not hold.
Therefore, $A\preceq^{\G}_{+} B$ does not hold.
\end{proof}

\subsection{Maximal elements and transitivity.}

\begin{proposition}\label{max_el_<}
Let $\G:M_{m\times n}\to\mathcal{P}(M_{n\times m})$ be a g-map and $A\in \Omega_\G$. Then $A$ is a maximal element in $(M_{m\times n},\le^{\G}_{+})$ if and only if $\operatorname{rank}(A)=\operatorname{min}\{m,n\}$.
\end{proposition}
\begin{proof}
Assume that $\operatorname{rank}(A)=\operatorname{min}\{m,n\}$ and $A\le^{\G}_{+} B$ for some $B\in M_{m\times n}$. Then $A\le^s B$ and there is $A^-\in \G(A)$ such that $A=AA^-BA^-A$. Since $A$ has maximal possible rank and $A\le^s B$, we conclude  that $\R(A)=\R(B)$ and $\N(A)=\N(B)$. It follows that $AA^-B=B$ and $B(I-A^-A)=0$, so $A=AA^-BA^-A=B$. Thus, $A$ is a maximal element in $(M_{m\times n},\le^{\G}_{+})$.

Conversely, assume that $\operatorname{rank}(A)<\operatorname{min}\{m,n\}$. Choose arbitrary $A^-\in\G(A)$ and set $P=AA^-$ and $Q=A^-A$. Since $\operatorname{rank}(A)<\operatorname{min}\{m,n\}$, $P\neq I_m$ and $Q\neq I_n$. Choose a nonzero vector $x\in \N(P)$ and a nonzero vector $y\in \N(Q^*)$. Let $B=A+xy^*$. Then $A\neq B$ and $A=PB=BQ=PBQ$, thus, $A\le^{\G}_{+} B$. Therefore, $A$ is not a maximal element in $(M_{m\times n},\le^{\G}_{+})$.
\end{proof}

In the following proposition we assume that the support $\Omega_\G$ contains all invertible matrices of order $n$; this is quite a natural assumption.

\begin{proposition}
Let $\G:M_{m\times n}\to\mathcal{P}(M_{n\times m})$ be a g-map and $A\in \Omega_\G$. If $\operatorname{rank}(A)=\operatorname{min}\{m,n\}$ then $A$ is a maximal element in ${(M_{m\times n},\preceq^{\G}_{+})}$.
Assume that $m=n$, $\G$ is semi-complete and all invertible matrices belong to the support $\Omega_\G$.
If $A$ is a maximal element in $(M_n,\preceq^{\G}_{+})$, then $A$ is invertible.
\end{proposition}
\begin{proof}
If $\operatorname{rank}(A)=\operatorname{min}\{m,n\}$ then $A$ is a maximal element in ${(M_{m\times n},\preceq^{\G}_{+})}$ by Proposition~\ref{max_el_<} and Proposition~\ref{preceq_sledi_<}.

Conversely, assume that $m=n$ and $A$ is not invertible. Fix $A^-\in\G(A)$. Since $\G$ is semi-complete, $G:=A^-AA^-\in \G_r(A)$. Set $P=AG$ and $Q=GA$. With respect to space decompositions $\mathbb{C}^n=\R(Q)\oplus\N(Q)$ and $\mathbb{C}^n=\R(P)\oplus\N(P)$, $A$ and $G$ have the following representations
$$A=\begin{bmatrix} A_1&0\\0&0\end{bmatrix}, \quad G=\begin{bmatrix} A_1^{-1}&0\\0&0\end{bmatrix},$$
where $A_1:\R(Q)\to \R(P)$ is invertible. Since $\operatorname{dim}\N(Q)=\operatorname{dim}\N(P)=n-\operatorname{rank}(A)$, we can find an invertible linear map $T:\N(Q)\to\N(P)$. Set $B=\begin{bmatrix} A_1&0\\0&T\end{bmatrix}$. Then $B$ is an invertible matrix and $A\neq B$.
By the assumption, $\G(B)=\{B^{-1}\}$, that is, $B^{-1}\in \G_r(B)$. Since $B^{-1}\{1\}=\{B\}$ and $GBG=G$, we conclude that $B^{-1}\{1\}\subseteq G\{1\}$. It follows that $G\le^- B^{-1}$. Obviously, $A\le^s B$, thus $A\preceq^{\G}_{+} B$.
Therefore, $A$ is not a maximal element in $(M_n,\preceq^{\G}_{+})$.
\end{proof}

\begin{proposition}\label{th_preceq_partial_order}
Let $\G:M_{m\times n}\to\mathcal{P}(M_{n\times m})$ be a g-map such that for $A\in\Omega_\G$, $\G(A)=\{A^-\}$ where $A^-\in A\{1,2\}$ is a fixed reflexive g-inverse of $A$. Then the relation $\preceq^{\G}_{+}$ is a partial order on $\Omega_\G$.
\end{proposition}
\begin{proof}
Since $A^-$ is reflexive, $\G$ is semi-complete, so by Proposition~\ref{reflanti}, the reflexivity and antisymmetry of $\preceq^{\G}_{+}$ hold.
The transitivity of $\preceq^{\G}_{+}$ follows directly from the definition of $\preceq^{\G}_{+}$ together with the transitivity of the relations $\le^s$ and $\le^-$.
\end{proof}

The situation is somewhat different regarding the relation $\le^{\G}_{+}$.
While the relations $\le^{\G}_{+}$ and $\preceq^{\G}_+$ do not coincide in general, Example~\ref{primer_G} demonstrates that they do coincide when $\mathcal{G}(A)=\{A^\dag\}$ or $\mathcal{G}(A)=\{A^\#\}$.
In the following theorem, we establish that they always coincide on the set of non-maximal elements provided that $\le^{\G}_{+}$ is transitive and $\mathcal{G}(A)$ consists of a fixed reflexive generalized inverse of $A$.

\begin{theorem}\label{th_<=preceq}
Let $\G:M_{m\times n}\to\mathcal{P}(M_{n\times m})$ be a g-map such that for $A\in\Omega_\G$, $\G(A)=\{A^-\}$ where $A^-\in A\{1,2\}$ is a fixed reflexive generalized inverse of $A$. Then the relation $\le^{\G}_{+}$ is transitive, i.e. $A\le^{\G}_{+} B$ and $B\le^{\G}_{+} C$ imply $A\le^{\G}_{+} C$ for all $A,B,C\in M_{m\times n}$,
if and only if the following implication holds for all $A,B\in\Omega_\G$:
\begin{equation}\label{<_sledi_preceq}
A\le^{\G}_{+} B \text{ and $B$ is not maximal in $(M_{m\times n},\le^{\G}_{+})$} \implies A\preceq^{\G}_{+} B.
\end{equation}
\end{theorem}
\begin{proof}
Assume that \eqref{<_sledi_preceq} holds for all $A,B\in\Omega_\G$ and let us show that $\le^{\G}_{+}$ is transitive. Suppose that $A\le^{\G}_{+} B$ and $B\le^{\G}_{+} C$. Then $A^-=A^-BA^-$ and $B^-=B^-CB^-$. If $C$ is not maximal in $(M_{m\times n},\le^{\G}_{+})$, then $B$ is not maximal either (since $B\le^{\G}_{+} C$), so $A\preceq^{\G}_{+} B$ and $B\preceq^{\G}_{+} C$. By Proposition~\ref{th_preceq_partial_order}, $A\preceq^{\G}_{+} C$, thus $A\le^{\G}_{+} C$, by Proposition~\ref{preceq_sledi_<}.
Suppose now that $C$ is a maximal element in $(M_{m\times n},\le^{\G}_{+})$. If $B=C$ then $A\le^{\G}_{+} C$. If $B\neq C$ then $B$ is not maximal, so $A\preceq^{\G}_{+} B$. Thus, $A^-=A^-BA^-=A^-BB^-=B^-BA^-$, by \eqref{preceq_karakterizacija} and
$$A^-CA^-=A^-BB^- C B^-BA^-=A^-BB^-BA^-=A^-BA^-=A^-.$$
Since $A\le^s C$, it follows that $A\le^{\G}_{+} C$.
Therefore, $\le^{\G}_{+}$ is a transitive relation.

Assume now that $\le^{\G}_{+}$ is transitive. Assume that $A\le^{\G}_{+} B$ where $B\in\Omega_\G$ is not maximal in $(M_{m\times n},\le^{\G}_{+})$. Then $A^-=A^-BA^-$, where $\G(A)=\{A^-\}$, and there exists $C\in M_{m\times n}$ such that $B\neq C$ and $B\le^{\G}_{+} C$. Let $\G(B)=\{B^-\}$ and define the following sets
$$V_A=\{T\in M_{m\times n} : A^-TA^-=0\}, \quad V_B=\{T\in M_{m\times n} : B^-TB^-=0\}$$ and
$$S=\{T\in M_{m\times n} : B^-TB^-=0,\, B\le^s B+T\}.$$
For $D\in M_{m\times n}$, since $B^-\in B\{1,2\}$, we see that $B\le^{\G}_{+} D$ if and only if $D-B\in S$.
Since $B\le^{\G}_{+} D$ implies $A\le^{\G}_{+} D$, while $A^-BA^-=A^-$, it follows that
$$S\subseteq V_A.$$
Also, $r:=\operatorname{rank}(B)<\operatorname{min}\{m,n\}$, by Proposition~\ref{max_el_<}.
The set $V_B$ is a subspace of $M_{m\times n}$ as it is equal to the kernel of a linear map $T\mapsto B^-TB^-$.
Let us show that $$\operatorname{span}(S)=V_B.$$
Since $B^-$ is a reflexive g-inverse of $B$, it follows, by Theorem~2.3.18 in \cite{knjiga} that there exist invertible matrices $K\in M_m$ and $L\in M_n$ such that
$$B=KD L^{-1} \quad \text{ and }\quad B^-=LD^- K^{-1},$$
where $D=\begin{bmatrix}I_r&0\\0&0\end{bmatrix}\in M_{m\times n}$,
$D^-=\begin{bmatrix}I_r&0\\0&0\end{bmatrix}\in M_{n\times m}$
and $I_r$ is the identity matrix of order $r$.
The map $\Phi(X)=K^{-1}XL$ is a linear bijection on $M_{m\times n}$ such that $\Phi(B)=D$, $X\le^s Y \,\Leftrightarrow\, \Phi(X)\le^s\Phi(Y)$ and $B^-TB^-=L\,D^-\Phi(T)D^-K^{-1}$. Hence $\Phi$ maps $V_B$ and $S$ onto the sets defined in the same way with $B$ and $B^-$ replaced by $D$ and $D^-$. Therefore, in order to prove that $\operatorname{span}(S)=V_B$, we may
assume that $K=I_m$ and $L=I_n$.
Note that $T\in V_B$ if and only if $T=\begin{bmatrix} 0&T_2\\T_3&T_4\end{bmatrix}$.
It follows that $\operatorname{dim}(V_B)=mn-r^2$ and a basis for the subspace
$V_B$ is
$$\mathfrak{B}=\{E_{i,j} : (i,j)\in\{1,\dots, m\}\times \{1,\dots, n\}\setminus \{1,\dots, r\}\times \{1,\dots, r\}\},$$
where $E_{i,j}$ denotes the matrix whose $(i,j)$-entry is $1$ and all other entries are $0$.
It is obvious that $E_{i,j}\in S$ when $r< i\leq m$ and $r< j\leq n$.
Assume the second case: $1\leq i\leq r$ and $r<j\leq n$. Fix $i_0$, $r<i_0\leq m$.
Let us show that $B\le^s B+E_{i,j}+E_{i_0,j}$.
For arbitrary $x\in\mathbb{C}^n$, let $\tilde{x}\in\mathbb{C}^n$ be defined by
$\tilde{x}_j=0$ and $\tilde{x}_k=x_k$ for $k\neq j$. Then $Bx=(B+E_{i,j}+E_{i_0,j})\tilde{x},$ and hence $\R(B)\subseteq\R(B+E_{i,j}+E_{i_0,j})$.
If $z\in \N(B+E_{i,j}+E_{i_0,j})$ then $E_{i_0,j}z=0$, so $z_j=0$. But then $E_{i,j}z=0$, and so $Bz=0$. Therefore, $\N(B+E_{i,j}+E_{i_0,j})\subseteq\N(B)$.
It follows that $E_{i,j}+E_{i_0,j}\in S$. Since $E_{i_0,j}\in S$, we conclude that $E_{i,j}\in \operatorname{span}(S)$. Assume the last case:
$r< i\leq m$ and $1\leq j\leq r$. Fix $j_0$, $r<j_0\leq n$. Let us show that $B\le^s B+E_{i,j}+E_{i,j_0}$.
It is obvious that $\N(B+E_{i,j}+E_{i,j_0})\subseteq\N(B)$.
For $x\in\mathbb{C}^n$, let $\tilde{x}\in\mathbb{C}^n$ be defined by $\tilde{x}_{j_0}=-x_j$ and $\tilde{x}_k=x_k$ for $k\neq j_0$. Then $Bx=(B+E_{i,j}+E_{i,j_0})\tilde{x}$. Therefore, $\R(B)\subseteq\R(B+E_{i,j}+E_{i,j_0})$.
It follows that $E_{i,j}+E_{i,j_0}\in S$. Since $E_{i,j_0}\in S$, we obtain that $E_{i,j}\in \operatorname{span}(S)$.
We have shown that
$$\mathfrak{B}\subseteq\operatorname{span}(S)\subseteq V_B.$$
Hence,  $\operatorname{span}(S)=V_B.$
Since $V_A$ is a subspace of $M_{m\times n}$ and $S\subseteq V_A$, we conclude that $V_B\subseteq V_A$, that is
\begin{equation}\label{pom_implikacija}
B^-TB^-=0 \implies A^-TA^-=0.
\end{equation}
Assume that $B^-GB^-=B^-$, $G\in M_{m\times n}$. Thus, $B^-(G-B)B^-=0$. By \eqref{pom_implikacija}, we obtain $A^-(G-B)A^-=0$. Since $A^-BA^-=A^-$, it follows that $A^-GA^-=A^-$. Therefore, $B^-\{1\}\subseteq A^-\{1\}$. From Theorem~\ref{minus_karakterizacije}, we conclude that $A^-\le^- B^-$, that is, $A\preceq^{\G}_{+} B$.
\end{proof}

\begin{remark}\label{remark_m=n}
If $m=n$ and $\Omega_\G$ contains all invertible matrices, then condition \eqref{<_sledi_preceq} in Theorem~\ref{th_<=preceq} is equivalent to
$$A\le^{\G}_{+} B \implies A\preceq^{\G}_{+} B \quad\text{for all } A,B\in\Omega_\G .$$
Indeed, by Proposition~\ref{max_el_<}, an element of $\Omega_\G$ is maximal in
$(M_n,\le^{\G}_{+})$ if and only if it is invertible.
Assume that $\G(A)=\{A^-\}\subseteq A\{1,2\}$, $A\le^{\G}_{+} B$ and $B$ is maximal in $(M_n,\le^{\G}_{+})$.
Then $A^-=A^-BA^-$ and $\G(B)=\{B^-\}$, where $B^-=B^{-1}$. Thus, $A^-=A^-BB^-=B^-BA^-$. By Theorem~\ref{minus_karakterizacije}, we conclude that $A^-\le^- B^-$, thus, $A\preceq^{\G}_{+} B$.
\end{remark}

\subsection{g-maps with non-singleton images}

When $\G(A)$ is not a singleton, the problem of characterizing the necessary and sufficient conditions under which $\le^{\G}_{+}$ and $\preceq^{\G}_{+}$ are partial orders becomes more involved.
In the following lemma, we provide a sufficient condition under which $\le^{\G}_{+}$ and $\preceq^{\G}_{+}$ coincide.

\begin{lemma}\label{th<=preceq}
Let $\G:M_{m\times n}\to\mathcal{P}(M_{n\times m})$ be a semi-complete g-map and $A,B\in \Omega_\G$. Suppose that $A\le^{\G}_{+} B$, that is, $A\le^s B$ and $A=AA^-BA^-A$ for some $A^-\in\G(A)$. If there exists $B^-\in\G(B)$ such that $B^-BA^-AA^-BB^-\in\G(A)$, then $A\preceq^{\G}_{+} B$.
\end{lemma}
\begin{proof}
Suppose that $A\le^s B$ and $A=AA^-BA^-A$ for some $A^-\in\G(A)$. By the assumption, there exists $B^-\in\G(B)$ such that $G:=B^-BA^-AA^-BB^-\in\G(A)$.
Since $A=BB^-A=AB^-B$, we have
\begin{align*}
GAG&=B^-BA^-AA^-(BB^- A B^-B)A^-AA^-BB^- \\
&= B^-BA^-AA^-AA^-AA^-BB^- = B^-BA^-AA^-BB^- = G,
\end{align*}
thus, $G\in \G_r(A)$.
Let $H:=B^-BB^-$. Since $\G$ is semi-complete, we have $H\in \G(B)\cap B\{2\}=\G_r(B)$. Also, $B\in H\{1\}$.
Furthermore,
\begin{align*}
GBG&=B^-BA^-AA^-BB^-\cdot B\cdot B^-BA^-AA^-BB^- \\
&=B^-BA^-(AA^-BA^-A)A^-BB^-=B^-BA^-AA^-BB^-=G
\end{align*}
and $GBH=GBB^-=G$, $HBG=B^-BG=G$.
By Theorem~\ref{minus_karakterizacije}, (i)$\Leftrightarrow$(vi) with $H^-=B$, it follows that $G\le^-H$, and thus, $A\preceq^{\G}_{+}B$.
\end{proof}

In the following theorem, we provide a sufficient condition under which $\le^{\G}_{+}$ and $\preceq^{\G}_{+}$ are partial orders.
\begin{theorem}\label{th_dovoljan_uslov}
Let $\G:M_{m\times n}\to\mathcal{P}(M_{n\times m})$ be a semi-complete g-map.
If the following condition holds for all $A,B\in \Omega_\G$, $A^-\in\G(A)$ and $B^-\in\G(B)$
\begin{equation}\label{dovoljan_uslov}
A\le^{\G}_{+} B \implies B^-BA^-AA^-BB^-\in\G(A),
\end{equation}
then $\le^{\G}_{+}$ and $\preceq^{\G}_{+}$ coincide and define a partial order on $\Omega_\G$.
\end{theorem}
\begin{proof}
Assume that \eqref{dovoljan_uslov} holds. By Lemma~\ref{th<=preceq} and Proposition~\ref{preceq_sledi_<}, we conclude that $\le^{\G}_{+}$ and $\preceq^{\G}_{+}$ coincide.
Let us show that $\le^{\G}_{+}$ and $\preceq^{\G}_{+}$ define a partial order on $\Omega_\G$.
By Proposition~\ref{reflanti}, it is enough to show that $\le^{\G}_{+}$ is transitive. Suppose that $A,B,C\in\Omega_\G$, $A\le^{\G}_{+} B$ and $B\le^{\G}_{+} C$. Thus $A\le^s B$, $B\le^s C$ and there exist $A^-\in \G(A)$ and $B^-\in \G(B)$ such that $A=AA^-BA^-A$ and $B=BB^-CB^-B$. We have $A\le^s C$. Let $G=B^-BA^-AA^-BB^-$. Then $G\in\G(A)$ and
\begin{align*}
GCG&=B^-BA^-AA^-(BB^-\cdot C\cdot B^-B)A^-AA^-BB^-\\
&=B^-BA^-(AA^-BA^-A)A^-BB^-
=B^-BA^-AA^-BB^-=G,
\end{align*}
so $A=AGA=AGCGA$. It follows that $A\le^{\G}_{+} C$.
\end{proof}

\begin{remark}\label{rem_not_necessary}
Condition \eqref{dovoljan_uslov} is not necessary for $\le^{\G}_{+}$ to be a partial order. Let $A=\operatorname{diag}(1,0,0)$, let
$$H_0=\begin{bmatrix}1&0&1\\0&0&0\\1&0&1\end{bmatrix}\in A\{1,2\},$$
and define $\G(A)=A\{1,2\}\setminus\{H_0\}$ and $\G(Z)=Z\{1,2\}$ for $Z\ne A$. Since $\G(Z)\subseteq Z\{1,2\}$ for every $Z\in M_3$, the g-map $\G$ is semi-complete.
We show that $\le^{\G}_{+}$ coincides with the plus order $\le^{+}$, and hence is a partial order.

Since $A\le^s B$ appears in both definitions, by Example~\ref{primer_G}(i) it suffices to show that $A=AH_0BH_0A$ implies $A=AXBXA$ for some $X\in A\{1,2\}$ with $X\ne H_0$.
For every $X\in A\{1,2\}$, we have $\operatorname{rank}(X)=\operatorname{rank}(A)=1$, so $X=pq^T$ for some $p,q\in\mathbb C^3$, and $AXA=A$ gives $p_1q_1=1$.
Since $pq^T=(\lambda p)(\lambda^{-1}q)^T$, we may assume that $p_1=q_1=1$. Conversely, every matrix $pq^T$ with $p_1=q_1=1$ belongs to $A\{1,2\}$. For such $X$ we have $AX=e_1q^T$ and $XA=pe_1^T$, where $e_1=\begin{bmatrix}1&0&0\end{bmatrix}^T$. Since $A=e_1e_1^T$, it follows that
$$AXBXA=e_1\,(q^TBp)\,e_1^T=(q^TBp)\,A.$$
Hence $A=AXBXA$ if and only if $q^TBp=1$. Let $B=[b_{ij}]$. For $X=H_0$ we have $p=q=\begin{bmatrix}1&0&1\end{bmatrix}^T$, so $A=AH_0BH_0A$ means that $b_{11}+b_{31}+b_{13}+b_{33}=1$.

Now let $q=\begin{bmatrix}1&0&1\end{bmatrix}^T$ and $p=\begin{bmatrix}1&p_2&p_3\end{bmatrix}^T$. Then
$$q^TBp=b_{11}+b_{31}+\alpha p_2+\beta p_3,\qquad \alpha=b_{12}+b_{32},\quad \beta=b_{13}+b_{33}.$$
If $(\alpha,\beta)\neq(0,0)$, put $p_2=\beta$ and $p_3=1-\alpha$. If $\alpha=\beta=0$, put $p_2=p_3=0$. In both cases $q^TBp=b_{11}+b_{31}+\beta=1$ and $(p_2,p_3)\neq(0,1)$. Therefore, $X=pq^T\in A\{1,2\}$, $X\neq H_0$ and $A=AXBXA$.

On the other hand, let $B=\operatorname{diag}(1,1,0)$, $A^-=A\in\G(A)$ and
$$B^-=\begin{bmatrix}1&0&1\\0&1&0\\1&0&1\end{bmatrix}\in B\{1,2\}=\G(B).$$
Then $A\le^s B$ and $AA^-BA^-A=A$, so $A\le^{\G}_{+}B$, but $B^-BA^-AA^-BB^-=H_0\notin\G(A)$. Thus \eqref{dovoljan_uslov} fails.
\end{remark}

\subsection{g-maps defined by equations (1)--(9)}

The \emph{core inverse} $A^{\core}$ and the \emph{dual core inverse} $A_{\core}$ of a matrix $A\in M_n$ were introduced by Baksalary and Trenkler~\cite{BT}. They exist if and only if $A\in M_n^\#$, in which case they are unique and
$$A^{\core}=A^\#AA^\dag, \qquad A_{\core}=A^\dag AA^\#.$$
For $A\in M_n$ consider the equations
\begin{align*}
(6)\;\;&XA^2=A, & (7)\;\;&AX^2=X,\\
(8)\;\;&A^2X=A, & (9)\;\;&X^2A=X.
\end{align*}
By~\cite{Core_Hilbert}, $X=A^{\core}$ if and only if $X$ satisfies (1), (2), (3), (6) and (7), and $X=A_{\core}$ if and only if $X$ satisfies (1), (2), (4), (8) and (9). In fact, equations (1), (3) and (7) already characterize $A^{\core}$ \cite[Theorem~2.1]{Wang_Liu}, and analogously, (1), (4) and (9) characterize $A_{\core}$.

The core and dual core inverses combine the Moore--Penrose and group inverses. In the same spirit, we introduce the following relations, which combine the diamond and square orders.

\begin{definition}\label{def_cir_relations}
Let $A,B\in M_n^\#$. We say that $A\le^\bigcirc B$ if $A\le^s B$ and $A^*A^2=A^*BA$, and that $A\le_\bigcirc B$ if $A\le^s B$ and
$A^2A^*=ABA^*$.
\end{definition}

We show in Corollary~\ref{posledica_core} that these relations are partial orders on $M_n^\#$.

Although the requirements in Theorem~\ref{th_dovoljan_uslov} may appear to be too restrictive, we now show that they are satisfied by all g-maps of the form $\G(A)=A\{1,S\}$, where $S\subseteq\{2,\dots,9\}$ is fixed and
$$A\{1,S\}=A\{1\}\cap\bigcap_{i\in S}A\{i\}.$$
Moreover, the condition $A\le^{\G}_{+}B$ in \eqref{dovoljan_uslov} can even be weakened to $A\le^s B$.
We first need the following lemma, which is probably known; parts (i) and (ii) can be found, e.g., in \cite[Chapter~2, Theorems~3 and~4]{BIG}. We include a short proof for completeness.

\begin{lemma}\label{lem:equations}
Let $A\in M_{m\times n}$ and $X\in A\{1\}$. Then
\begin{enumerate}[label=(\roman*)]
\item $X\in A\{3\}$ if and only if $AX=AA^\dagger$;
\item $X\in A\{4\}$ if and only if $XA=A^\dagger A$.
\end{enumerate}
If $m=n$ and $A\{1,i\}\neq\emptyset$ for some $i\in\{5,6,7,8,9\}$, then
$A\in M_n^\#$, and
\begin{enumerate}[label=(\roman*),resume]
\item $X\in A\{5\}$ if and only if $AX=XA=AA^\#$;
\item $X\in A\{6\}$ if and only if $XA=AA^\#$;
\item $X\in A\{7\}$ if and only if $X\in A\{2\}$ and $XA=AA^\#$;
\item $X\in A\{8\}$ if and only if $AX=AA^\#$;
\item $X\in A\{9\}$ if and only if $X\in A\{2\}$ and $AX=AA^\#$.
\end{enumerate}
\end{lemma}
\begin{proof}
Throughout, $X\in A\{1\}$.

(i) If $(AX)^*=AX$, then $AX=AA^\dag AX=(AXAA^\dag)^*=AA^\dag$. The converse is clear. The proof of (ii) is analogous to the proof of (i).

If $X$ satisfies equations (5), (6), (7), (8) or (9), then $A=A^2X$, $A=XA^2$, $A=A^2X^2A$, $A=A^2X$ or $A=AX^2A^2$, respectively. In each case $\operatorname{rank}(A^2)=\operatorname{rank}(A)$, so $A\in M_n^\#$.

(iv) If $XA^2=A$, then $XA=XA^2A^\#=AA^\#$. Conversely, $XA=AA^\#$ gives $XA^2=AA^\#A=A$.
The proof of (vi) is analogous to the proof of (iv).

(iii) If $AX=XA$, then $A=A^2X=XA^2$, so by (iv) and (vi), $AX=XA=AA^\#$. The converse is clear.

(v) If $AX^2=X$, then $\operatorname{rank}(X)=\operatorname{rank}(X^2)$, so $X\in M_n^\#$. Thus $XAX=XAX^2X^\#=X^2X^\#=X$ and
$XA=AX^2A=A^\#A(AX^2)A=A^\#AXA=A^\#A$.
Conversely, let $X\in A\{2\}$ and $XA=AA^\#$. Then $X=XAX=AA^\#X$, and
$AX^2=AX(XA)X=(AXA)A^\#X=AA^\#X=X$.
(vii) is proved analogously to (v).
\end{proof}

For $A\in M_{m\times n}$ (and $A\in M_n^\#$ whenever the group inverse is involved) put
\begin{align*}
&\PP_\dagger(A)=\{AA^\dagger\},\quad \PP_\#(A)=\{AA^\#\},\quad \PP_+(A)=\PP(A),
\\
&\QQ_\dagger(A)=\{A^\dagger A\},\quad \QQ_\#(A)=\{AA^\#\},\quad \QQ_+(A)=\QQ(A),
\end{align*}
and, for $x,y\in\{\dagger,\#,+\}$, define
$$A\le^{(x,y)}B \iff A\le^s B \text{ and } A=PBQ
\text{ for some } P\in\PP_x(A),\ Q\in\QQ_y(A).$$

\begin{theorem}\label{teorema_klasifikacija}
Let $S\subseteq\{2,\dots,9\}$ and let $\G(A)=A\{1,S\}$, where $m=n$ if
$S\cap\{5,\dots,9\}\neq\emptyset$. Define $x,y\in\{\dag,\#,+\}$ by
$$x=\begin{cases}
\dag, & 3\in S,\\
\#, & 3\notin S,\ S\cap\{5,8,9\}\neq\emptyset,\\
+, & \text{otherwise},
\end{cases}$$
and
$$y=\begin{cases}
\dag, & 4\in S,\\
\#, & 4\notin S,\ S\cap\{5,6,7\}\neq\emptyset,\\
+, & \text{otherwise}.
\end{cases}$$
Then the following hold:
\begin{itemize}
\item[(i)] $\Pi_\G(A)=\PP_x(A)\times\QQ_y(A)$ for every $A\in\Omega_\G$.
If $3\in S$ and $S\cap\{5,8,9\}\neq\emptyset$, or $4\in S$ and
$S\cap\{5,6,7\}\neq\emptyset$, then $\Omega_\G$ consists of EP matrices.
\item[(ii)] $\G$ is semi-complete, and for all $A,B\in\Omega_\G$,
$A^-\in\G(A)$ and $B^-\in\G(B)$,
\begin{equation}\label{dovoljan_uslov_oslabljeni}
A\le^s B \implies B^-BA^-AA^-BB^-\in\G(A).
\end{equation}
\item[(iii)] The relations $\le^{\G}_{+}$, $\preceq^{\G}_{+}$ and $\le^{(x,y)}$ coincide on $\Omega_\G$, and they are partial orders on $\Omega_\G$. In particular, if $\Omega_\G$ consists of EP matrices, then $\le^{\G}_{+}$ is the restriction of $\le^{(x,y)}$ to the set of EP matrices.
\end{itemize}
\end{theorem}

\begin{proof}
(i) Let $A\in\Omega_\G$. If $3\in S$ and $S\cap\{5,8,9\}\neq\emptyset$,
then for $X\in\G(A)$ Lemma~\ref{lem:equations} gives
$AA^\dag=AX=AA^\#$, so $A$ is EP; the case $4\in S$,
$S\cap\{5,6,7\}\neq\emptyset$ is similar. Conversely, every EP matrix
belongs to $\Omega_\G$, since $A^\dag=A^\#$ satisfies equations (1)--(9).

Let $X\in A\{1\}$. By Lemma~\ref{lem:equations}, each equation $(j)$,
$j\in S$, is equivalent to one or more of the conditions $AX=AA^\dag$
($j=3$), $AX=AA^\#$ ($j\in\{5,8,9\}$), $XA=A^\dag A$ ($j=4$),
$XA=AA^\#$ ($j\in\{5,6,7\}$) and $X\in A\{2\}$ ($j\in\{2,7,9\}$). By the
choice of $x$ and $y$, and since $AA^\dag=A^\dag A=AA^\#$ in the cases
where $A$ is EP, we obtain
\begin{equation}\label{karakterizacija}
X\in\G(A)\iff
\begin{aligned}[t]
&AX\in\PP_x(A),\ XA\in\QQ_y(A),\\
&\text{and } X\in A\{2\} \text{ if } S\cap\{2,7,9\}\neq\emptyset.
\end{aligned}
\end{equation}
In particular, $\Pi_\G(A)\subseteq\PP_x(A)\times\QQ_y(A)$ and
\begin{equation}\label{refl}
X\in A\{1,2\},\ AX\in\PP_x(A),\ XA\in\QQ_y(A)\implies X\in\G(A).
\end{equation}
Conversely, let $(P,Q)\in\PP_x(A)\times\QQ_y(A)$. By
Lemma~\ref{remark_P,Q}, $X=QA^-P$, where $A^-\in A\{1\}$ is arbitrary, is
a reflexive inverse of $A$ with $AX=P$ and $XA=Q$, so $X\in\G(A)$ by
\eqref{refl} and $(P,Q)\in\Pi_\G(A)$.

(ii) Let $A,B\in\Omega_\G$, $A^-\in\G(A)$ and $B^-\in\G(B)$. The matrix $H=A^-AA^-$ is a reflexive inverse of $A$ with $AH=AA^-$ and $HA=A^-A$, so $H\in\G(A)$ by (i) and \eqref{refl}. Thus $\G$ is semi-complete.

Assume now that $A\le^s B$ and let $G=B^-BA^-AA^-BB^-$. Since $A=BB^-A=AB^-B$, we have $AGA=AA^-A=A$,
$$AG=AA^-BB^-,\qquad GA=B^-BA^-A,$$
and $GAG=G$ as in the proof of Lemma~\ref{th<=preceq}, so $G\in A\{1,2\}$. If $x=\dag$, then $AA^-=AA^\dag$ and $BB^-=BB^\dag$ by Lemma~\ref{lem:equations}, so $AG=AA^\dag BB^\dag=(BB^\dag AA^\dag)^*=(AA^\dag)^*=AA^\dag$. If $x=\#$, then $AA^-=AA^\#$ and $BB^-=BB^\#$, so $AG=A^\#ABB^\#=A^\#A=AA^\#$, since $A=AB^\#B$. Hence $AG\in\PP_x(A)$, and similarly $GA\in\QQ_y(A)$. By \eqref{refl}, $G\in\G(A)$, i.e. \eqref{dovoljan_uslov_oslabljeni} holds.

(iii) By (ii), condition \eqref{dovoljan_uslov} holds, so by Theorem~\ref{th_dovoljan_uslov}, ${\le^{\G}_{+}}={\preceq^{\G}_{+}}$ is a partial order on $\Omega_\G$. By Proposition~\ref{proposition_pairs}(i) and (i), ${\le^{\G}_{+}}={\le^{(x,y)}}$ on $\Omega_\G$.
\end{proof}

\begin{remark}\label{napomena_24}
By \eqref{karakterizacija}, for $X\in A\{1\}$ the condition $X\in A\{1,S\}$ depends only on the condition imposed on $AX$ (none,
$AX=AA^\dag$ or $AX=AA^\#$), on the condition imposed on $XA$ (none, $XA=A^\dag A$ or $XA=AA^\#$), and on whether $X\in A\{2\}$ is required. If two different conditions are imposed on $AX$ or on $XA$, then $\Omega_\G$ consists of EP matrices, on which $AA^\dag=A^\dag A=AA^\#$; in this case only three possibilities remain ($AX$ restricted, $XA$
restricted, or both). Hence, for $m,n\ge 2$, one easily checks that there are exactly $(3\cdot 3+3)\cdot 2=24$ different g-maps of the form $\G(A)=A\{1,S\}$; the $8$ of them with $S\subseteq\{2,3,4\}$ are defined for all $m,n$, while the remaining $16$ require $m=n$.
For instance, $A\{1,5\}=A\{1,6,8\}$, $A\{1,2,5\}=A\{1,7,8\}=\{A^\#\}$ and $A\{1,2,3,6,7\}=A\{1,3,7\}=\{A^{\core}\}$, which recovers \cite[Theorem~2.1]{Wang_Liu}.
\end{remark}

Thus, although there are $24$ different g-maps $A\{1,S\}$, by
Theorem~\ref{teorema_klasifikacija} they induce only nine relations,
namely $\le^{(x,y)}$, $x,y\in\{\dag,\#,+\}$. If $m\neq n$, only the four
relations with $x,y\in\{\dag,+\}$ occur. In the cases described in the
theorem where $\Omega_\G$ consists of EP matrices, we obtain restrictions
of these relations to EP matrices. Of course, other g-maps lead to further
relations; for instance, $\G(A)=\{A^\dag_{M,N}\}$ given by the weighted
Moore--Penrose inverse, or g-maps in which different sets $S$ are used for
different matrices. These nine relations are listed in
Table~\ref{tabela_relacija}. Each condition in the table is obtained from
$A=PBQ$ by multiplying it on the left by $A^*$ if $P=AA^\dag$, or by $A$
if $P=AA^\#$, and on the right by $A^*$ if $Q=A^\dag A$, or by $A$ if
$Q=AA^\#$; multiplying by $(A^*)^\dag$ or $A^\#$ on the corresponding
sides gives back $A=PBQ$.

Five of these relations have already appeared: the diamond order ($x=y=\dag$) \cite{Baksalary_Hauke_diamond}, the plus order ($x=y=+$) \cite{Arias_Maestripieri}, and the square order and the relations $\le^\bigcirc$, $\le_\bigcirc$ introduced in Section~2 and Definition~\ref{def_cir_relations}.
For instance, $S=\{3,7\}$ gives $A\{1,S\}=\{A^{\core}\}$ and $\le^{(\dag,\#)}=\le^\bigcirc$, while the classes $A\{1,2,3\}$ and
$A\{1,2,4\}$ of 1MP and MP1 inverses~\cite{1MP} induce $\le^{(\dag,+)}$ and $\le^{(+,\dag)}$, respectively.
Since $\PP_x(A)\subseteq\PP(A)$ and $\QQ_y(A)\subseteq\QQ(A)$,
$$A\le^{(x,y)}B \implies \bigl(A\le^{(x,+)}B \text{ and } A\le^{(+,y)}B\bigr)\implies A\le^+B.$$

\begin{table}[ht]
\centering
\renewcommand{\arraystretch}{1.3}
\begin{tabular}{c|c|c|c}
 & $y=\dag$ & $y=\#$ & $y=+$ \\ \hline
$x=\dag$ & $\le^\diamondsuit$ & $\le^\bigcirc$ & $\le^{(\dag,+)}$ \\
 & $AA^*A=AB^*A$ & $A^*A^2=A^*BA$ & $A^*A=A^*BQ$ \\ \hline
$x=\#$ & $\le_\bigcirc$ & $\le^\square$ & $\le^{(\#,+)}$ \\
 & $A^2A^*=ABA^*$ & $A^3=ABA$ & $A^2=ABQ$ \\ \hline
$x=+$ & $\le^{(+,\dag)}$ & $\le^{(+,\#)}$ & $\le^+$ \\
 & $AA^*=PBA^*$ & $A^2=PBA$ & $A=PBQ$
\end{tabular}
\medskip
\caption{The relations $\le^{(x,y)}$, $x,y\in\{\dag,\#,+\}$. In each cell,
$A\le^{(x,y)}B$ if and only if $A\le^s B$ and the displayed equation holds
for some $P\in\PP(A)$ and $Q\in\QQ(A)$; if $x=\#$ or $y=\#$, then
$A\in M_n^\#$.}
\label{tabela_relacija}
\end{table}

\begin{corollary}\label{posledica_core}
For $A,B\in M_n^\#$ the following assertions are equivalent:
\begin{itemize}
\item[(i)] $A\le^\bigcirc B$, that is, $A\le^s B$ and $A^*A^2=A^*BA$;
\item[(ii)] $A\le^s B$ and $A^{\core}=A^{\core}BA^{\core}$;
\item[(iii)] $A\le^s B$ and $A^{\core}\le^-B^{\core}$.
\end{itemize}
Analogously, $A\le_\bigcirc B$ if and only if $A\le^s B$ and $A_{\core}=A_{\core}BA_{\core}$, if and only if $A\le^s B$ and $A_{\core}\le^-B_{\core}$. In particular, $\le^\bigcirc$ and $\le_\bigcirc$ are partial orders on $M_n^\#$.
\end{corollary}

\begin{proof}
Let $\G(A)=A\{1,3,7\}=\{A^\core\}$, i.e., $S=\{3,7\}$. Then $\Omega_\G=M_n^\#$, $x=\dag$ and $y=\#$. By Table~\ref{tabela_relacija}, (i) is $A\le^{(\dag,\#)}B$; by the characterization after Definition~\ref{def<}, (ii) is $A\le^{\G}_{+}B$; and since $\G(A)=\G_r(A)$, (iii) is $A\preceq^{\G}_{+}B$. By Theorem~\ref{teorema_klasifikacija}(iii), these three relations coincide and are a partial order on $M_n^\#$. For $\le_\bigcirc$ take $S=\{4,9\}$, i.e., $\G(A)=\{A_\core\}$.
\end{proof}

\begin{remark}
Note that the condition $A^\core \le^- B^\core$ implies $\R(A)\subseteq\R(B)$, but does not imply $\N(B)\subseteq\N(A)$.
Analogously, the condition $A_\core \le^- B_\core$ implies $\N(B)\subseteq\N(A)$, but does not imply $\R(A)\subseteq\R(B)$.

Indeed, if $A^\core \le^- B^\core$ then $A^\core=B^\core BA^\core$. Therefore
$$A=A^\core A^2=B^\core BA^\core A^2=B^\core B A=BB^\# A,$$
thus, $\R(A)\subseteq\R(B)$.

For a counterexample to $\N(B)\subseteq\N(A)$, let $A=\begin{bmatrix}1&1\\0&0 \end{bmatrix}$ and $B=\begin{bmatrix}1&0\\0&0 \end{bmatrix}$. Then $A^\core=B^\core=B$, so $A^\core \le^- B^\core$. But it is evident that $\N(B)\not\subseteq\N(A)$.

The second statement can be proved similarly.
\end{remark}

\subsection{Comparison with Mitra's g-based orders}

By Proposition~\ref{proposition_pairs}, for a semi-complete g-map $\G$, Mitra's relation $\le^{\G}$ and the relation $\le^{\G}_{+}$ are determined by the same set of pairs $\Pi_\G(A)$:
$$A\le^{\G}B \iff A=PB=BQ, \qquad
A\le^{\G}_{+}B \iff A\le^s B \text{ and } A=PBQ,$$
for some $(P,Q)\in\Pi_\G(A)$. For $x,y\in\{\dag,\#,+\}$ we define
$$A\le^{[x,y]}B \iff A=PB=BQ \text{ for some } P\in\PP_x(A),\ Q\in\QQ_y(A).$$

\begin{corollary}\label{posledica_Mitra}
Let $S$, $\G$, $x$ and $y$ be as in Theorem~\ref{teorema_klasifikacija}.
Then $\le^{\G}$ coincides with $\le^{[x,y]}$ on $\Omega_\G$.
\end{corollary}
\begin{proof}
By Theorem~\ref{teorema_klasifikacija}(ii), $\G$ is semi-complete.
The statement follows from Proposition~\ref{proposition_pairs}(ii) and
$\Pi_\G(A)=\PP_x(A)\times\QQ_y(A)$ (Theorem~\ref{teorema_klasifikacija}).
\end{proof}
The relations $\le^{[x,y]}$ are well-known partial orders (see \cite{knjiga} and, for the core and dual core orders, \cite{BT}; see also Table~\ref{tabela_Mitra}). In particular, the left and right star orders were described by 1MP and MP1 inverses in \cite[Theorem~5]{Rakic1MP}. Corollary~\ref{posledica_Mitra} shows that these nine orders are the only relations induced in Mitra's sense by the g-maps $A\{1,S\}$.

The conditions in Table~\ref{tabela_Mitra} are obtained as follows. Multiplying $A=AA^\dagger B$ and $A=AA^\# B$ on the left by $A^*$ and $A$, and $A=BA^\dagger A$ and $A=BAA^\#$ on the right by $A^*$ and $A$, gives $A^*A=A^*B$, $A^2=AB$, $AA^*=BA^*$ and $A^2=BA$, respectively; multiplying by $(A^\dag)^*$ or $A^\#$ on the corresponding sides gives back the original equations. Moreover, if $A=PB$ for some $P\in\mathsf P(A)$, then $A=BQ$ for some $Q\in\mathsf Q(A)$ if and only if $\R(A)\subseteq\R(B)$.
Indeed, if $\R(A)\subseteq\R(B)$, then $A=BX$ for some $X\in M_n$. Let $A^-\in A\{1\}$ and $Q=XA^-A$. We have $AX=PBX=PA=A$, so $Q^2=XA^-AXA^-A=Q$ and $AQ=A$. Since also $Q=XA^-A$, it follows that $\N(Q)=\N(A)$, so $Q\in\mathsf Q(A)$, and $BQ=BXA^-A=AA^-A=A$.
Similarly, if $A=BQ$ for some $Q\in\mathsf Q(A)$, then $A=PB$ for some $P\in\mathsf P(A)$ if and only if $\N(B)\subseteq\N(A)$.

\begin{table}[ht]
\centering
\renewcommand{\arraystretch}{1.3}
\begin{tabular}{c|c|c|c}
 & $y=\dag$ & $y=\#$ & $y=+$ \\ \hline
$x=\dag$ & $\le^*$ & $\le^{\core}$ & $\le^*_\ell$ \\
 & $A^*A=A^*B$ & $A^*A=A^*B$ & $A^*A=A^*B$ \\
 & $AA^*=BA^*$ & $A^2=BA$ & $\mathcal R(A)\subseteq\mathcal R(B)$ \\ \hline
$x=\#$ & $\le_{\core}$ & $\le^\#$ & $\le^\#_\ell$ \\
 & $A^2=AB$ & $A^2=AB$ & $A^2=AB$ \\
 & $AA^*=BA^*$ & $A^2=BA$ & $\mathcal R(A)\subseteq\mathcal R(B)$ \\ \hline
$x=+$ & $\le^*_r$ & $\le^\#_r$ & $\le^-$ \\
 & $AA^*=BA^*$ & $A^2=BA$ & $A=PB=BQ$ \\
 & $\mathcal N(B)\subseteq\mathcal N(A)$ & $\mathcal N(B)\subseteq\mathcal N(A)$ &
\end{tabular}
\medskip
\caption{The relations $\le^{[x,y]}$, $x,y\in\{\dag,\#,+\}$: the star,
core, left star, dual core, sharp, left sharp, right star, right sharp and
minus orders. In each cell, $A\le^{[x,y]}B$ if and only if the displayed
conditions hold (for some $P\in\PP(A)$ and $Q\in\QQ(A)$ in the last cell);
if $x=\#$ or $y=\#$, then $A\in M_n^\#$.}
\label{tabela_Mitra}
\end{table}

Tables~\ref{tabela_relacija} and~\ref{tabela_Mitra} have the same structure, and by Proposition~\ref{<_ne_sledi_preceq}(i), $A\le^{[x,y]}B$ implies $A\le^{(x,y)}B$; for instance, the star, sharp and minus orders imply the diamond, square and plus orders, respectively. There is, however, an essential difference. In Mitra's theory, $P$ and $Q$ appear in separate equations, so that
${\le^{[x,y]}}={\le^{[x,+]}}\cap{\le^{[+,y]}}$; e.g., the star order is the intersection of the left and right star orders. In our setting $P$ and $Q$ are coupled in $A=PBQ$, and only ${\le^{(x,y)}}\subseteq{\le^{(x,+)}}\cap{\le^{(+,y)}}$ holds. The inclusion
is strict for $x,y\in\{\dag,\#\}$: if
$$A=\begin{bmatrix}1&0\\0&0\end{bmatrix},\quad
B=\begin{bmatrix}0&1\\1&0\end{bmatrix},\quad
P=\begin{bmatrix}1&1\\0&0\end{bmatrix},\quad
Q=\begin{bmatrix}1&0\\1&0\end{bmatrix},$$
then $A=A^2=A^*$, $P\in\PP(A)$, $Q\in\QQ(A)$ and $A=ABQ=PBA$, so $A$ is below $B$ in each of the relations $\le^{(x,+)}$ and $\le^{(+,y)}$, $x,y\in\{\dag,\#\}$, but $ABA=0\neq A$, so $A\not\le^{(x,y)}B$.

Several results of this section are counterparts of results of Mitra. Condition~\eqref{dovoljan_uslov} in Theorem~\ref{th_dovoljan_uslov} plays the role of Mitra's (T)-condition: $A\le^{\G}B$ implies $\{B^-AB^-: B^-\in\G(B)\}\subseteq\G(A)$ \cite[Theorem~2.1]{Unified_theory}. For g-maps whose images are singletons consisting of reflexive g-inverses, Mitra \cite[Theorem~2.6]{Unified_theory} proved that $\le^{\G}$ is a partial order if and only if
$$A\le^{\G}B \text{ and } B \text{ is not maximal} \implies B^-AB^-=A^-,$$
which should be compared with Theorem~\ref{th_<=preceq}, where the conclusion is $A^-\le^-B^-$.
As in Remark~\ref{remark_m=n}, for $m=n$ the assumption that $B$ is not maximal can be omitted (see \cite[Corollary~7.2.29]{knjiga}).
Finally, while \cite[Theorem~7.5.5]{knjiga} provides a class of g-maps inducing partial orders, Theorem~\ref{teorema_klasifikacija} and Corollary~\ref{posledica_Mitra} determine all relations induced by the g-maps $A\{1,S\}$ in both theories.
Mitra's theory was extended in \cite{Rakic_unified} to the so-called FD rings, which generalize the ring of linear operators on a finite-dimensional vector space.

\section{Conclusion}

In this paper we introduced the square partial order $\le^\square$ on $M_n^\#$, the group inverse analogue of the diamond order, characterized by $A^\#\le^-B^\#$. We described it by block representations, related it to the reverse order law for the group inverse, and showed that, unlike the Moore--Penrose case, the conditions $A\le^-B$ and $(B-A)^\#=B^\#-A^\#$ do not characterize the sharp order, but lead to the new partial orders $\le^\boxplus_\ell$, $\le^\boxplus_r$ and $\le^\boxplus$.

We also developed a framework, parallel to Mitra's theory of g-based orders, in which the plus order plays the role of the minus order. For a g-map $\G$, both Mitra's relation and the relation $\le^{\G}_{+}$ are determined by the pairs of idempotents $(AA^-,A^-A)$, $A^-\in\G(A)$; the difference between them is the difference between the conditions $A=PB=BQ$ and $A=PBQ$. For the g-maps $A\{1,S\}$, both theories give nine relations. In Mitra's theory these are well-known orders, whereas in our setting, besides the diamond and plus orders, we obtain the square order, the orders $\le^\bigcirc$ and $\le_\bigcirc$, and four further relations.

It is natural to ask which of these results can be extended beyond matrices. The first such setting is that of bounded linear operators on Hilbert spaces, to which the minus order was originally extended in \cite{AntezanaCorachStojanoff_Minus,Semrl_Minus}; another is that of rings, where many results of this paper can be extended by the standard techniques used, e.g., for the sharp and core orders in \cite{Rakic_core_annihilators,Rakic_star_sharp_core}. In both settings some care is needed; for instance, the proof of Theorem~\ref{additivity_group_th} relies on the Cayley--Hamilton theorem, and the proof of Theorem~\ref{th_<=preceq} essentially uses linear algebra. It would be interesting to know which of these results remain valid for operators and in rings, where the plus and diamond orders have already been studied \cite{Arias_Maestripieri,CirulisDiamond,LebtahiPatricioThome,MarovtMosic_plus}.

Several further questions remain open.
First, necessary and sufficient conditions under which $\le^{\G}_{+}$ and $\preceq^{\G}_{+}$ are partial orders are not known for g-maps whose images are not singletons; Theorem~\ref{th_dovoljan_uslov} only provides a sufficient condition, which is not necessary by Remark~\ref{rem_not_necessary}. Second, the orders $\le^\bigcirc$, $\le_\bigcirc$ and the relations $\le^{(x,+)}$, $\le^{(+,y)}$ deserve a more detailed study. Finally, the square order is defined only for group invertible matrices. Since a relation defined only in terms of the Drazin inverse does not see the nilpotent part of a matrix, a natural extension to matrices of arbitrary index should also take the nilpotent part into account, as in the C-N partial order \cite{Mitra_Hartwig}, see also \cite[Section~4.4]{knjiga}.

\section*{Acknowledgements}
This research was financially supported by the Ministry of Science, Technological Development and Innovations of the Republic of Serbia (Contract No. 451-03-34/2026-03/200109).

\section*{ORCID}

\noindent Dragan S. Raki\'c - \url{https://orcid.org/0000-0002-3556-2885}


\begin{thebibliography}{99}

\bibitem{AntezanaCorachStojanoff_Minus}
J. Antezana, G. Corach, D. Stojanoff,
Bilateral shorted operators and parallel sums,
\textit{Linear Algebra Appl.}
\textbf{414} (2006) 570--588.

\bibitem{Arias_Maestripieri}
M. L. Arias, A. Maestripieri,
On partial orders of operators,
\textit{Ann. Funct. Anal.}
\textbf{14} (2023) 21.

\bibitem{Baksalary_aditivity_MP}
J. K. Baksalary,
A relationship between the star and minus orderings,
\textit{Linear Algebra Appl.}
\textbf{82} (1986) 163--167.

\bibitem{Baksalary_Hauke_diamond}
J. K. Baksalary, J. Hauke,
A further algebraic version of Cochran’s theorem and matrix partial orderings,
\textit{Linear Algebra Appl.}
\textbf{127} (1990) 157--169.

\bibitem{BT}
O. M. Baksalary and G. Trenkler,
Core inverse of matrices,
\textit{Linear Multilinear Algebra}
\textbf{58}(6) (2010) 681--697.

\bibitem{BIG}
A. Ben-Israel, T. N. E. Greville,
Generalized Inverses: Theory and Applications,
2nd Edition, Springer (2003).

\bibitem{BurgosMGMC}
M. Burgos, A. C. M\'arquez-Garc\'ia, A. Morales-Campoy,
Maps preserving the diamond partial order,
\textit{Appl. Math. Comput.}
\textbf{296} (2017) 137--147.

\bibitem{CirulisDiamond}
J. C\={\i}rulis,
The diamond partial order for strong Rickart rings,
\textit{Linear Multilinear Algebra},
\textbf{65}(1) (2017) 192--203.

\bibitem{Drazin_natural_structures}
M. P. Drazin,
Natural structures on semigroups with involution,
\textit{Bull. Amer. Math. Soc.}
\textbf{84}(1) (1978) 139--141.

\bibitem{FerreyraLevisMaharanaOrquera}
D. E. Ferreyra, F. E. Levis, G. Maharana, V. Orquera,
New characterizations of the diamond partial order,
\textit{J. Comput. Appl. Math.}
\textbf{476} (2026) 117087.

\bibitem{FerreyraMalikCore}
D. E. Ferreyra, S. B. Malik,
Some new results on the core partial order,
\textit{Linear Multilinear Algebra},
\textbf{70}(18) (2022) 3449--3465.

\bibitem{Fongi_Gonzalez}
G. Fongi, M. C. Gonzalez,
Moore-Penrose inverse and partial orders on Hilbert space operators,
\textit{Linear Algebra Appl.}
\textbf{674} (2023) 1--20.

\bibitem{Ha}
R. E. Hartwig,
How to partially order regular elements,
\textit{Math. Japon.}
\textbf{25} (1980) 1--13.

\bibitem{Hartwig_Styan}
R. E. Hartwig, G. P.H. Styan,
On some characterizations of the “star” partial ordering for matrices and rank subtractivity,
\textit{Linear Algebra Appl.}
\textbf{82} (1986)  145--161.

\bibitem{1MP}
M. V. Hern\'andez, M. B. Lattanzi, N. Thome,
From projectors to 1MP and MP1 generalized inverses and their induced partial orders,
\textit{RACSAM}
\textbf{115} (2021) 148.

\bibitem{HernandezLattanziThome_Diamond}
M. V. Hern\'andez, M. B. Lattanzi, N. Thome,
On diamond partial order, one-sided star partial orders, and 1MP-inverses,
\textit{RACSAM}
\textbf{118} (2024) 145.

\bibitem{LebtahiPatricioThome}
L. Lebtahi, P. Patr\' icio, N. Thome,
The diamond partial order in rings,
\textit{Linear Multilinear Algebra},
\textbf{62}(3) (2014) 386--395.

\bibitem{MalikRuedaThome}
S. B. Malik, L. Rueda, N. Thome,
Further properties on the core partial order and other matrix partial orders,
\textit{Linear Multilinear Algebra},
\textbf{62}(12) (2014) 1629--1648.

\bibitem{MarovtMihelic}
J. Marovt, K. Miheli\v c,
On sets of elements in Rickart rings induced by partial orders,
\textit{Appl. Math. Comput.}
\textbf{315} (2017) 555--563.

\bibitem{MarovtMosic_plus}
J. Marovt, D. Mosi\'c, I. Cremer,
On some generalized inverses and partial orders in $*$-rings,
\textit{J. Algebra Appl.}
\textbf{22}(12) (2023) 2350256.

\bibitem{Unified_theory}
S. K. Mitra,
Matrix partial orders through generalized inverses: unified theory,
\textit{Linear Algebra Appl.}
\textbf{148} (1991) 237--263.

\bibitem{Mitra_sharp_order}
S. K. Mitra,
On group inverses and the sharp order,
\textit{Linear Algebra Appl.}
\textbf{92} (1987) 17--37.

\bibitem{knjiga}
S. K. Mitra, P. Bhimasankaram, S.B. Malik,
\textit{Matrix partial orders, shorted operators and applications},
World Scientific (2010).

\bibitem{Mitra_Hartwig}
S. K. Mitra, R. E. Hartwig,
Partial orders based on outer inverses,
\textit{Linear Algebra Appl.} \textbf{176} (1992) 3--20.

\bibitem{Nambooripad}
K. S. S. Nambooripad,
The natural partial order on a regular semigroup,
\textit{Proc. Edinburgh Math. Soc.}
\textbf{23} (1980) 249--260.

\bibitem{Rakic_core_annihilators}
D. S. Raki\'c,
Generalization of sharp and core partial order using annihilators,
\textit{Banach J. Math. Anal.}
\textbf{9}(3) (2015) 228--242.

\bibitem{Rakic_star_sharp_core}
D. S. Raki\'c, D. S. Djordjevi\'c,
Star, sharp, core and dual core partial order in rings with involution,
\textit{Appl. Math. Comput.}
\textbf{259} (2015) 800--818.

\bibitem{Rakic_unified}
D. S. Raki\'c, D. S. Djordjevi\'c,
Partial orders in rings based on generalized inverses -- unified theory,
\textit{Linear Algebra Appl.}
\textbf{471} (2015) 203--223.

\bibitem{Core_Hilbert}
D. S. Raki\'c, N. \v C. Din\v ci\'c, D. S. Djordjevi\'c,
Core inverse and core partial order of Hilbert space operators,
\textit{Appl. Math. Comput.}
\textbf{244} (2014) 283--302.

\bibitem{Rakic1MP}
D. S. Raki\'c, M. Z. Ljubenovi\'c,
1MP and MP1 inverses and one-sided star orders in a ring with involution,
\textit{RACSAM}
\textbf{117} (2023) 13.

\bibitem{Semrl_Minus}
P. \v Semrl,
Automorphisms of $B(H)$ with respect to minus partial order,
\textit{J. Math. Anal. Appl.}
\textbf{369} (2010) 205--213.

\bibitem{Wang_Liu}
H. Wang, X. Liu,
Characterizations of the core inverse and the core partial ordering,
\textit{Linear Multilinear Algebra}
\textbf{63}(9) (2015) 1829--1836.

\bibitem{YangJi}
L. Yang, G. Ji,
Bounds for the diamond partial order in $B(H)$,
\textit{Linear Multilinear Algebra}
\textbf{69}(8) (2021) 1415--1421.


\end{thebibliography}
\end{document}